
\let\he=\heads
\let\coo=\comments
\def\hh #1\par{\underbar{\he#1}\vskip1pc}
\def\\{\noindent}
\def\qed{\hskip 12.75cm $\heartsuit$\medbreak}
\def\gu{[GuSh]}
\def\vv{\par\rm}
\def\pt #1. {\medbreak\noindent{\bf#1. }\enspace\sl}
\def\pd #1. {\medbreak \noindent{\bf#1. }\enspace}
\def\dd{Definition\ } 
\def\sh{[Sh]}
\def\bb #1{{\bar #1}}
\def\fo{\ ||\kern -4pt-}
\def\la{\lambda}
\def\de{\delta}
\def\ale{\aleph_0}
\def\sc{semi--club}
\def\sb{\subseteq}
\def\sbb{\subset}
\def\st{such that\ }
\def\om{\omega}
\def\al{\alpha}
\def\be{\beta}
\def\ga{\gamma}
\def\imp{\Rightarrow}
\def\ind{induction hypothesis\ }
\def\lan{\langle}
\def\ran{\rangle}
\def\magbas{\magnification 1200 \baselineskip 1.2pc}
\def\sec{\vskip 2pc}
\def\fps{$n$-formally possible set of theories}
\def\fpss{$n$-formally possible sets of theories}
\def\tk{$T_k$}
\def\fff #1{\bb #1^{\thinspace\frown}}
\def\ff #1{#1^{\thinspace\frown}}
\def\all{\forall}
\def\ex{\exists}
\def\qmu{$Q$\llap{\lower5pt\hbox{$\sim$}\kern1pt}$_\mu$}
\font\title=cmbx12
\font\by=cmr8
\font\author=cmr10
\font\adress=cmsl10
\font\abstract=cmr8
\def\tit #1\par{\centerline{\title #1}}
\def\bby{\centerline{\by BY}\par}
\def\aut #1\par{\centerline{\author #1}}
\def\abs{\abstract\centerline{ABSTRACT}\par}
\def\abss #1\par{\abstract\midinsert\narrower\narrower\noindent #1\endinsert}
\sec
\magbas
\tit PEANO ARITHMETIC MAY NOT BE INTERPRETABLE \par
\tit IN THE MONADIC THEORY OF ORDER \par 
\vskip.5in
\bby
\aut SHMUEL LIFSCHES and SAHARON SHELAH\footnote*{The second author would like 
to thank the U.S.--Israel Binational Science Foundation for 
partially supporting this research. Publ. 471}
\par
{\adress Institute of Mathematics, The Hebrew University of Jerusalem, 
Jerusalem, Israel} \par
\vskip.5in
\abs
\abss Gurevich and Shelah have shown that Peano Arithmetic cannot be 
interpreted in the monadic second-order theory of short chains 
(hence, in the monadic second-order theory of the real line). 
We will show here that it is consistent that there is no 
interpretation even in the monadic second-order theory of all chains. \par
\vskip.6in
%
%
%
%
\hh 0. Introduction \par \rm
A {\sl reduction} 
of a theory $T$ to a theory $T^*$ is an algorithm, associating a sentence 
$\varphi^*$ in the language of $T^*$, to each sentence $\varphi$ in the 
language of $T$, in such a way that:\ \  $T \vdash \varphi\ \  $if and only if
$\ \ T^* \vdash \varphi^*$.

\\Although reduction is a powerful method of proving undecidability results,
it lacks in establishing any semantic relation between the theories.

\\A {\sl (semantic) interpretation} of a theory \ $T$ in a theory \ $T'$ is a 
special case of reduction in which models of \ $T$ are defined inside models
of \ $T'$. \par
\\It is known (via reduction) that the monadic theory of order and the monadic 
theory of the real line are complicated at least as Peano Arithmetic, (In \sh\  
this was proven from ZFC+MA and in [GuSh1] from ZFC), and even as second order 
logic ([GuSh2], [Sh1], for the monadic theory of order). 
Moreover, second order logic was shown to be interpretable in the monadic 
theory of order ([GuSh3]) but this was done by using a weaker, 
non--standard form of interpretation: into a Boolean valued model. Using 
standard interpretation ([GMS]) it was shown that it is consistent that the 
second--order theory of $\om_2$ is interpretable in the monadic theory 
of $\om_2$. \par
\\On the other hand, by \gu, Peano Arithmetic is not interpretable 
in the monadic theory of short chains, and in particular in the monadic 
theory of the real line.  \par
\\More details and Historical background can be found in [Gu]. \par
The previous results leave a gap concerning the question whether it is  
provable from ZFC that Peano Arithmetic is interpretable in the monadic 
theory of order. 
\ In this paper we fill the gap and show that the previous results are the 
best possible, by proving:
\proclaim Theorem. There is a forcing notion $P$ such that in $V^P$,
Peano Arithmetic (in fact a much weaker theory) is not interpretable and even 
not weakly interpretable in the monadic second-order theory of chains. 
\vv \par
\\From another point of view the theorem may be construed as presenting the 
strength of the interpretation method by showing that although 
Peano Arithmetic is recursive in the monadic theory of order, it is not  
interpretable in it.  \medbreak
\\In the proof we use notations and definitions from \sh\ and \gu\ but 
although we omit some proofs, it is self contained. \par
\\We start by defining in $\S1$ the notion of interpretation.  Althogh this 
notion is not uniform in the literature, our notion of weak interpretation 
seems to follow from every reasonable definition. \par
\\In $\S2$ we define partial theories and present the relevant results about 
them from \sh. \par
\\In $\S3$ we define a theory $T$, easily interpretable in Peano Arithmetic, 
with the following axioms:

(a) $\all x\ex y\all z[p(z,y)\leftrightarrow z=x]$ \par
(b) $\all x\all y\ex u\all z[p(z,u)\leftrightarrow(p(z,x)\vee p(z,y))]$ \par
(c) $\ex x\all y[\neg p(y,x)]$ \vv
\\Assuming there is a chain $C$ that interprets $T$, we show that the 
interpretation `concentrates' on an initial segment $D\sb C$. \medbreak
\\The main idea in the proof is that of shuffling subsets $X,Y\sb C$: 
Given a partition of $C$, \ $\lan S_j: j\in J \ran$ and a subset $a\sb J$, 
the shuffling of $X$ and $Y$ with respect to $J$ and $a$ is the set: \ 
$\bigcup_{j\in a}(X\cap S_j)\cup\bigcup_{j\not\in a}(Y\cap S_j)$. \ We show in 
$\S4$ and $\S5$ that under suitable circumstances (in particular, 
if $a$ is a `\sc'), partial theories are preserved under shufflings. \par
\\We use a simple class forcing $P$, defined in $\S5$, to obtain a universe  
$V^P$ in which generic semi--clubs are added to every suitable partition. \par
\\The contradiction to the assumption that an interpretation exists in 
$V^P$ can be roughly described as follows: \ 
We start with an interpreting chain $C$. 
The interpretation defines an equivalence relation between subsets of $C$,
and we choose a large enough number of nonequivalent subsets.  
We fix a partition of $C$ and after some manipulations we are left with 3 
ordered pairs of nonequivalent subsets of $C$.
We shuffle each pair $U, V$ with respect to a generic \sc\  $a$, added by the 
forcing, and get a new subset 
which is equivalent to $U$. (This uses the preservation of partial theories 
undershufflings). 
But a condition $p\in P$ that forces these 
equivalences determines only a bounded subset of $a$. We show that we could 
have got the same results if we had shuffled the pairs with respect to the 
complement of $a$. Thus for each pair $U, V$, $p$ forces that the `inverse' 
shuffling is also equivalent to $U$. We conclude by showing that one of 
the shufflings is equivalent to $V$ as well, and get a contradiction since 
$U$ and $V$ were not equivalent.
\sec 
%
%
%
%
\hh 1. The notion of interpretation \par
The notion of semantic interpretation of a theory $T$ in a theory $T'$ is not 
uniform. Usually it means that models of $T$ are defined inside models of 
$T'$ but the definitions vary with context. Here we will define
the notion of interpretation of one first order theory in another following 
the definitions and notatins of \gu.
%
%
%
%
\pd Remark. The idea of our definition is that in every model of $T'$ (or 
maybe of some extension $T''$ if $T'$ is not complete) we can define a model 
of $T$. An alternative definition could demand that {\sl every} model of 
$T$ is interpretable in a model of $T'$ (As in [BaSh]). Actually we need a 
weaker notion than the one we define and this seems to follow from every 
reasonable definition of semantic interpretation. \ We will show that it is 
consistent that no chain $C$ interprets Peano arithmetic. We even allow 
parameters from $C$ in the interpreting formulas. \ Thus, our notion is:   
``A model of $T'$ defines (with parameters) a model of $T$\ ''. 
We call this notion ``{\sl Weak Interpretation''} \vv
%
%
%
%
\pd \dd 1.1. Let $\sigma$ be a signature $\lan{\rm P_1,P_2},\ldots\ran$ where 
each $\rm P_i$ is a predicate symbol of some arity $\rm r_i$, in the language
$L=L(\sigma)$. 
An {\sl interpretation} of $\sigma$ in a first order language $L'$ is a 
sequence  \par
\\$ I = \lan d, U(\bb v_1,\bb u), E(\bb v_1,\bb v_2,\bb u),
P'_1(\bb v_1,\ldots \bb v_{r_1},\bb u),
P'_2(\bb v_1,\ldots \bb v_{r_2},\bb u), \ldots \ran $ \ \ where: \par
(a) $d$ is a positive integer (the {\sl dimension}); \par
(b) $U(\bb v_1,\bb u)$ and $E(\bb v_1,\bb v_2,\bb u)$ are $L'$-formulas 
    (the {\sl universe {\rm and the} equality formulas});   \par
(c) each $P'_i(\bb v_1,\ldots \bb v_{r_i},\bb u)$ is an $L'$-formula (the 
    {\sl interpretation} of $P_i$); \par
(d) $\bb v_1, \bb v_2\ldots$ are disjoint $d$-tuples of distinct variables of 
$L'$; \par 
(e) $\bb u$ is a finite sequence (standing for the parameters of the 
    interpretation). \vv
%
%
%
%
\pd \dd 1.2. Let $\sigma$, $L'$ and $I$ be as in 1.1. Fix a function 
that associates each $L$ variable $v$ with a $d$-tuple $v'$ of distinct $L'$ 
variables in such a way that if $u$ and $v$ are different $L$-variables 
then the tuples $u'$ and $v'$ are disjoint. \par
\\ We define, by induction, the {\sl $I$-translation} 
$\varphi'$  of an arbitrary $L$-formula $\varphi$:           \par
(a) $(x=y)' = E(x',y')$.  \par
(b) If P is a predicate symbol of arity $r$ in $L$, then 
    ${\rm P}(x_1\ldots x_r)' = {\rm P}'(x'_1\ldots x'_r)$.   \par
(c) $(\neg\varphi)' = \neg(\varphi')$, and 
    $(\varphi\wedge\psi)' =  (\varphi'\wedge\psi')$.         \par
(d) $(\forall x)\varphi(x)' = (\forall x')\lbrack U(x')\rightarrow\varphi'(x')
    \rbrack$, and $(\exists x)\varphi(x)' = (\exists x')\lbrack 
    U(x')\wedge\varphi'(x')\rbrack$.                         \vv
%
%
%
%
\pd \dd 1.3. Let $T$ and $T'$ be first order theories such that the 
signature of $T$ consists of predicate symbols, and $T'$ is consistent
and complete. Let $I$ be an interpretation of the signature of $T$ in 
$L(T')$, and let $U(x)$ be the universe formula of $I$.          \par
\\$I$ is an {\sl interpretation} of $T$ in $T'$ if:                \par 
(a) the formula $\exists xU(x)$ is a theorem of \ $T'$, and      \par
(b) the $I$ translation of every closed theorem of $T$ is a 
    theorem of $T'$.                                             \par
\\$T$ is {\sl interpretable} in $T'$ if there is an interpretation of $T$ in 
$T'$. \vv
%
%
%
%
\pd \dd 1.4. Let $T, T'$ and $U(x)$ as in 1.3. except  $T'$ may be 
incomplete. Let $T''$ be the extension of $T'$ by an additional axiom 
$\exists xU(x)$.     \par
\\$I$ is an {\sl interpretation} of $T$ in $T'$ if:       \par
(a) $T''$ is consistent, and                            \par
(b) the $I$ translation of every closed theorem of $T$ 
    is a theorem of $T''$.                              \par
\\$T$ is {\sl interpretable} in $T'$ if there is an interpretation of $T$ in 
$T'$. \vv
%
%
%
%
\pd Remark 1.5. 1) The definitions are easily generalized to the case that 
$\sigma(T)$ consists also of function symbols, see \gu. \par
\\2) Definitions 1.3 and 1.4 make sense in case there are no parameters
in the interpretation. \vv
%
%
%
%
\pd \dd 1.6. Let $\sim$ be an equivalence relation on a non empty set $A$, 
and let $R$ be a relation of some arity $r$ on $A$. We say that \enspace
$\sim$ {\sl respects} $R$ \enspace if for all elements 
$a_1,\ldots,a_r, \ b_1,\ldots,b_r$ of $A$,  \par
$\lbrack R(a_1,\ldots,a_r)\ \&\ (a_1\sim b_1)\ \&\ \ldots\ \&\ (a_r\sim b_r)
\rbrack$ \ implies \ $R(b_1,\ldots,b_r)$. \vv
%
%
%
%
\pd \dd 1.7. Let $\sigma$, $I$ and $L'$ be as in def. 1.1. Let $M$ be a model   
for $L'$ and  \par
(a) $U^* = \{x:x$ is a $d$-tuple of elements of $M$ and $U(x)$ 
    holds in $M\}$;                     \par
(b) $E^* = \{(x,y):\ x,y\in U^*$ and $E(x,y)$ holds in $M\}$; \ and   \par
(c) if P is a predicate symbol of arity $r$ in $\sigma$, then            \par       
    \ \ \ \ ${\rm P}^* = \{(x_1,\ldots,x_r):$ each $x_i$ belongs to $U^*$ and 
    ${\rm P}'((x_1,\ldots,x_r)$ holds in $M\}$.                          \par
\\The interpretation $I$ {\sl respects} the structure \ $M$ if \ $U^*$ is 
not empty, \ $E^*$ is an equivalence relation, and  \ $E^*$ respects every \ 
${\rm P}^*$. \ (The definition is easily generalized when we allow parameters 
in $I$). \vv
%
%
%
%
\pt Lemma 1.8. Any interpretation of a first-order theory $T$ in a consistent  
complete first order theory $T'$ respects every model of $T$. \vv
\pd Proof. Easy                                              
\vv \qed 
%
%
%
%
\pd \dd 1.9. Let $\sigma$, $I$ and $L'$ be as in \dd 1.1. and let $M, U^*, 
E^*$ and ${\rm P}^*$ be as in \dd 1.7. \ We suppose that $I$ respects $M$ 
and define a Model for $L$ which will be called the {\sl I-image} of $M$ and 
will be denoted $I(M)$.                                                 \par
\\Elements of $I(M)$ are equivalence classes 
$x/E^* = \{y\in U^*: xE^*y\}$ of $E^*$ (where $x$ ranges over 
$U^*$). If P is a predicate symbol of arity $r$ in $\sigma$ then P is 
interpreted in $I(M)$ as the relation $\{(x_1/E^*,\ldots, x_r/E^*):
(x_1,\ldots,x_r)\in {\rm P}^*\}$. Again, we may allow parameters in $I$ and                                        
slightly modify this definition. \vv
%
%
%
%
\pt Lemma 1.10. Let $I = (d, U(v_1), E(v_1,v_2),\ldots)$ be an interpretation
of a signature $\sigma$ in the first order language of a structure $M$. 
Suppose that $I$ respects $M$. \ Let:      \par
$\varphi(v_1,\ldots,v_l)$ be an arbitrary $L(\sigma)$-formula, 
$\varphi'(v'_1,\ldots,v'_l)$ its $I$-translation,                \par
$U^* = \{x:x$ is a $d$-tuple of elements of $M$ and $U(x)$ holds in $M\}$,\par
$E^* = \{(x,y): x,y\in U^*$ and $E(x,y)$ holds in $M\}$,   \par
$x_1,\ldots,x_l$ belong to $U^*$.      \par
\\Then, \ $\varphi'(x_1,\ldots,x_l)$ holds in $M$ if and only if 
$\varphi(x_1/E^*,\ldots,x_l/E^*)$  holds in $I(M)$.                        \vv
\pd Proof. By induction on $\varphi$. 
\vv \qed
%
%
%
%
\\So we can conclude:
\pt Theorem 1.11. If $I$ is an interpretation of a first-order theory $T$ in 
the first-order theory of a structure $M$, {\sl then} the $I$-image of $M$ 
is a model for $T$.          \vv
\pd Proof. Let $\varphi$ be any closed theorem of $T$. Since $I$ interpretes
$T$ in the theory of $M$, the $I$-translation $\varphi'$ of $\varphi$ holds
in $M$. By Lemma 1.10, $\varphi$ holds in $I(M)$. 
\vv \qed
%
%
%
%
\pd Remark 1.12. The notion of interpretation presents a connection between 
theories: It implies that models of a theory $T$ are defined inside models 
of the interpretating theory $T'$. (Assuming $T'\vdash(\exists x)U(x)$,   
for every $M\models T$, $I(M)$ is a model of $T$).  But rephrasing a previous 
remark we demand less: In our world $V^P$ we will show that there is no
model $M$ of (actually a weaker theory than) Peano Arithmetic, and no
chain $C$ (= a model of the monadic theory of order), and an interpretation 
$I$, \st the $I$-image of $C$ is isomorphic to $M$. This will hold even if we
allow parameters in the interpreting formulas in $I$. This leads to the 
following definition:
%
%
%
%
\pd \dd 1.13. $T'$ {\sl weakly interprets} $T$ if there is a model $M$ of 
$T'$ and an interpretation $I$ of the signature of $T$, respecting $M$, 
maybe with parameters from $M$ appearing in $I$, \st 
$I(M)$ is a model of $T$.  \medbreak
\\From now on, whenever we write `interpretation' we will mean weak 
interpretation  in the sense of the previous definition.
\sec 
%
%
%
%
\hh 2. Partial Theories \par
In this section we will define 3 kinds of partial theories following \sh: \ 
$Th^n$ (definition 2.3) which is the theory of formulas with monadic 
quantifier depth n, \ $ATh^n$ (definition 2.11) which is the n-theory of 
segments (and by 2.10 `many' segments have the same theory), and \ $WTh^n$ 
which gives information about stationary subsets of the chain. 
The last two theories are naturally defined for well ordered chains only, 
but by embeding a club in the chain we can modify them so that they 
can be applied also to general chains. \par
\\The main result of this section  states roughly that for every n there is 
an m \st $WTh^m$ and $ATh^m$ determine $Th^n$ (theorem 2.15). \par
%
%
%
%
\pd \dd 2.1. The {\it monadic second-order theory} of a chain $C$ is the 
theory of $C$ in the language of order enriched by adding variables for sets 
of elements, atomic formulas of the form ``$x\in Y$''and the quantifier 
$(\exists Y)$ ranging over subsets. Call this language $L$. \vv
\pd Remark. We can identify the monadic theory of $C$ with the 
first order theory of the associated structure \ 
$$ C' = \lan {\cal P}(C), \sbb, <, \emptyset \ran $$   
where ${\cal P}(C)$ is the power set of $C$, and $<$ is the binary relation 
$\{(\{x\},\{y\}): x,y$ are elements of $C$ and $x<y$ in $C \}$.  \vv
%
%
%
%
\pd Notation. The universe of a model $M$ will be denoted $|M|$. Let $x,y,z$
be individual variables; $X,Y,Z$ set variables; $a,b,c$ elements; $A,B,C$ 
sets. Bar denotes a finite sequence, like $\bb a$, and $l({\bar a})$ it's  
length. We write e.g. ${\bar a}\in M$ and $\bb A\sb M$ instead of 
${\bar a}\in |M|^{l({\bar a})}$, or  $\bb A\in {\cal P}(|M|)^{l(\bb A)}$   \vv
%
%
%
%
\pd \dd 2.2. For any $L$-model $M$, $\bb A\in {\cal P}(M),\ \bb a\in|M|$, 
and a natural number $n$ define
$$t = th^n(M,\bb A,\bb a)$$ 
by induction on $n$:

\\for $n=0$:  \ \ 
$t = \{\varphi(X_{l_1},\ldots,x_{j_1},\ldots):\varphi(X_{l_1},\ldots,x_{j_1},
\ldots)$ is an atomic formula in $L$, $M\models 
\varphi[A_{l_1},\ldots,a_{j_1},\ldots]\}$. \par
\\for $n=m+1$: \ \ 
$t = \{th^m(M,\bb A,\bb a^\wedge b):\bb b\in |M|\}$. \vv  
%
%
%
%
\pd \dd 2.3. For any $L$-model $M$, $\bb A\in {\cal P}(M)$, and a natural 
number $n$ define
$$T = Th^n(M,\bb A)$$ 
by induction on $n$:

\\for $n=0$:\ \ \ \ \ \ \ \ \ $T = th^2(M,\bb A)$. \par
\\for $n=m+1$:
\ \ \ $T = \{Th^m(M,\bb A^\wedge B):B\in {\cal P}(M)\}$. \medbreak
\pd Remark. By $Th^0(M,\bb A)$ we can tell which subset is a singleton, 
so we can proceed to quantify only over subsets. \vv
%
%
%
%
\pt Lemma 2.4. (A) For every formula $\psi(\bb X)\in L$ there is an $n$  
such that from $Th^n(M,\bb A)$ we can find effectively whether 
$M \models \psi(\bb X)$. \par
\\(B) For every $n$ and $m$ there is a set $\Psi = 
\{\psi_l(\bb X):l<l_0(<\omega), l(\bb X)=m\} \subset L$ such that for any 
$L$-models $M, N$ and $\bb A\in {\cal P}(M)^m,\bb B\in {\cal P}(N)^m$ the 
following hold: \par
(1) $Th^n(N,\bb B)$ can be computed from $\{l<l_0:N\models 
    \psi_l[\bb B]\}$ \par
(2) $Th^n(M,\bb A)=Th^n(N,\bb B)$ if and only if for any
    $l<l_0$, \ \ $M\models \psi_l[\bb A]\leftrightarrow N\models 
    \psi_l[\bb B]$. \vv
\pd Proof. \rm In \sh, \ Lemma 2.1  (Note that our language $L$ is finite).
\vv \qed
%
%
%
%
\pt Lemma 2.5. For given $n,m$, \ each $Th^n(M,\bb A)$  
is hereditarily finite, (where $l(\bb A)=m, \ M$\ is an $L$-model), \ and we 
can effectively compute the set of formally possible $Th^n(M,\bb A)$. \vv
\pd Proof. In \sh, \ Lemma 2.2 
\vv \qed  
%
%
%
%
\pd \dd 2.6. If $C,D$ are chains then $C+D$ \ is any chain that can be split
into an initial segment isomorphic to $C$ and a final segment isomorphic to 
$D$.  \par
\\If $\langle C_i:i<\al\rangle$ is a sequence of chains then 
$\sum_{i<\al}C_i$ \ is any chain $D$ that is the concatenation of segments  
$D_i$, such that each $D_i$ \ is isomorphic to $C_i$. \vv
%
%
%
%
\pt Theorem 2.7 (composition theorem). \par
(1) If \ $l(\bb A)=l(\bb B)=l$, and
$$Th^m(C,\bb A) = Th^m(C',\bb A')$$ 
and 
$$Th^m(D,\bb B) = Th^m(D',\bb B')$$
then 
$$Th^m(C+D,A_0\cup B_0,\ldots,A_{l-1}\cup B_{l-1}) = 
Th^m(C'+D',A'_0\cup B'_0,\ldots,A'_{l-1}\cup B'_{l-1}).$$ 
\\(2) If \ $Th^m(C_i,\bb {A_i}) = Th^m(D_i,\bb {B_i})$ for each $i<\al$, then
$$Th^m\Bigl( \sum_{i<\al}C_i,\ \cup_i A_{1,i},\ldots,\cup_i A_{l-1,i}\Bigr) = 
Th^m\Bigl( \sum_{i<\al}D_i,\ \cup_i B_{1,i},\ldots,\cup_i B_{l-1,i}\Bigr).$$ 
\vv
\pd Proof. By \sh\ Theorem 2.4 (where a more general theorem is proved), 
or directly by induction on $m$. 
\vv \qed
%
%
%
%
\pd Notation 2.8. (1) $Th^m(C,A_0,\ldots,A_{l-1}) + Th^m(D,B_0,\ldots,
B_{l-1})$ \ \ is \ $Th^m(C+D,A_0\cup B_0,\ldots,A_{l-1}\cup B_{l-1})$. \par
\\(2) $\sum_{i<\al}Th^m(C_i,\bb {A_i})$ \ \ is  \ 
  $Th^m(\sum_{i<\al}C_i,\ \cup_{i<\al}A_{1,i},\ldots,\cup_{i<\al} 
  A_{l-1,i})$. \par
\\(3) If $D$ is a subchain of $C$ and $X_1,\ldots,X_{l-1}$ \ are subsets of 
  $C$ then $Th^m(D,X_0,\ldots,X_{l-1})$ abbreviates 
  $Th^m(D,X_0\cap D,\ldots,X_{l-1}\cap D)$.  \vv
\vskip 25pt
The following definitions and results apply to well ordered chains 
(i.e. ordinals), later we will modify them.
\vskip 20pt
%
%
%
%
\pd \dd 2.9. For $a\in (M,\bb A)$ let 
$$th(a,\bb A) = \{x\in X_i:a\in A_i\}\cup\{x\not\in X_i:a\not\in A_i\}.$$ 
So it is a finite set of formulas. \medbreak
For $\al$ \ an ordinal with $cf(\al)>\omega$, \ let $D_\al$ denote
the filter generated by the closed unbounded subsets of $\al$.
%
%
%
%
\pt Lemma 2.10. If the cofinality of $\al$ is $>\omega$, \ then for every 
$\bb A\in {\cal P}(\al)^m$ \ there is a closed unbounded subset $J$ \ of $\al$
such that: for each $\beta<\al$,\ all the models 
$$\{(\al, \bb A)|_{[\beta,\gamma)}:\gamma\in J,\ cf(\gamma)=\omega,\  
\gamma>\beta\}$$
have the same monadic theory. \par
\pd Proof. In \sh \ Lemma 4.1. 
\vv \qed
%
%
%
%
\pd \dd 2.11. $ATh^n(\beta,\ (\al,\bb A))$ for\ $\beta<\al,
\al$ a limit ordinal of cofinality $>\omega$ \ is\ 
$Th^n(\al,\bb A)|_{[\beta,\gamma)})$\ for every $\gamma\in J,\ 
\gamma>\beta,\ cf(\gamma)=\omega$;\ Where\ $J$\  is from Lemma 2.10.  \vv
\pd Remark. As\ $D_\al$\ is a filter, the definition does not depend on  
the choice of \ $J$.  \vv
%
%
%
%
\pd \dd 2.12. We define \ \ $WTh^n(\al, \bb A)$: \par
\\(1) if\ $\al$\ is a successor or has cofinality \ $\omega$,\  it is 
\ $\emptyset$;  \par
\\(2) otherwise we define it by induction on \ $n$:     \par
for \ $n=0$:\ \  
$WTh^0(\al, \bb A) = \bigl\{t:\{\beta<\al: th(\beta,\bb A)=t\}$ 
\ is a stationary subset of \ $\al\ \bigr\}$;           \par
for \ $n+1$:\ \ 
$WTh^{n+1}(\al, \bb A) = \{\langle S^{\bb A}_1(B),S^{\bb A}_2(B)\rangle
:B\in {\cal P}(\al)\}$         
\par
\\Where:    \par
$S^{\bb A}_1(B) = WTh^n(\al, \bb A,B)$,    \par
$S^{\bb A}_2(\bb B) = \bigl\{\langle t,s\rangle:\{\beta<\al:   
WTh^n((\al, \bb A,B)|_\beta)=t,\                
\ th(\beta, \bb A^\wedge B)=s\} $ 
\ is a stationary subset of \ $\al\bigr\} $.  \vv
\pd Remark. Clearly, if we replace \ $(\al,\bb A)$\ by a submodel whose 
universe is a club subset of \ $\al$,\ $WTh^n(\al, \bb A)$\ 
will not change. \vv
%
%
%
%
\pd \dd 2.13. Let \ $cf(\al)>\omega,\ M = (\al,\bb A)$\  and we define
the model \ $g^n(M) = \bigr(\al,g^n(\bb A)\bigl)$. \par
\\Let:\ \ $\bigl(g^n(\bb A)\bigr)_s = \{\beta<\al:s=ATh^n
(\beta,(\al,{\bb A}))\}$     \par
\\let $m=l(\bb A)$\ and \ $T(n,m)$\ := the set of formally possible 
\ $Th^n(M,\bb B)$, where $l(\bb B)=m$.  We define:
$$g^n(\bb A) := \lan\ldots, \bigr(g^n(\bb A)\bigl)_s,\ldots
\ran_{s\in T(n,m)}.$$  \vv
%
%
%
%
\pt Lemma 2.14. (A) \ $g^n(\al,\bb A)$ \ is a partition of $\al$. \par
\\(B) \ $g^n(\al,\bb A^\wedge \bb B)$ is a refinement of 
$g^n(\al,\bb A)$ and we can effectively correlate the parts. \par
\\(C) \ $g^{n+1}(\al,\bb A)$ is a refinement of   
$g^n(\al,\bb A)$ and we can effectively correlate the parts. \vv
\pd Proof. Easy. 
\vv \qed
\vskip 20pt
The next theorem shows that the (partial) monadic theories can be 
computed from \ $ATh$\ and \ $WTh$\ and is the main tool for showing that 
the monadic theories are preserved under shufflings of subsets. \par 
%
%
%
%
\pt Theorem 2.15. If \ $cf(\al)>\omega$,\  then for each \ $n$\ there
is an\ $m=m(n)$\ such that if:  \par
\\$t_1=WTh^m\bigl(\al,g^m(\al,\bb A)\bigr),
\ t_2=ATh^m\bigl(0,(\al,\bb A)\bigr)$  \par
then we can effectively compute \ $Th^n(\al,\bb A)$\ from\ $t_1,t_2$. \vv
\pd Proof. By \sh, Thm. 4.4. 
\vv \qed
%
%
%
%
\pd Notation 2.16. We will denote \ $ \langle t_1,t_2 \rangle $\ 
from Thm. 2.15 by \ $WA^{m(n)}$. \vv  \sec
%
%
%
%
In \sh\ the partial theories $ATh$ and $WTh$ were defined only to 
well ordered chains. We will show now how we can modify our definitions and 
apply them to general chains of cofinality $>\om$. 
The only loss of generality is that we assume that we can find in every chain 
$C$ a closed cofinal sequence. This does not hurt us because if a chain $C$ 
interprets a theory $T$, then there is a chain $C^c$ that interprets $T$, with 
this property and all we have to pay is maybe adding an additional parameter 
to the interpreting formulas.
The proofs of the results are easy generalizations of the original proofs.
\pd Notation 2.17. Let $C$ be a chain of cofinality $\la>\om$, and 
$J^*=\lan \be_i : i<\la \ran$ be a closed cofinal subchain of $C$. 
Fix a club subset of $\la$, \ $J=\lan \al_i : i<\la \ran$ \st $\al_0=0$ 
and for simplicity $cf(\al_{i+1})=\om$, and let $h\colon J^*\rightarrow J$ 
be an isomorphism, $h(\be_i)=\al_i$.
For a fixed $n$ and $\bb A\sb{\cal P}(C)^d$, denote by $s_i$ the theory 
$Th^n\big((C,\bb A)|_{[\be_i,\be_{i+1})}\big)$. Using these notations we can 
generalize the definitions and facts concerning $ATh$ and $WTh$: \vv
%
%
%
%
\pt Lemma 2.10*. If the cofinality of $C$ is $>\omega$, \ then for every 
$\bb A\in {\cal P}(C)^d$ \ there is a subchain $J^{**}\sb J^*$ 
\st $h''(J^{**})=J'\sb J$ is a club subset of $\la$, with $0\in J'$, 
and \st for each $i<\la$,\ all the models 
$$\{ (C,\bb A)|_{[\be_i,\be_j)} : j>i, \be_j\in J^{**}, cf(h(\be_j))=\om \}$$
have the same monadic theory. \vv     \qed
\pd Remark. We could have chosen $J$ to be all $\la$. The definitions 
and the results do not depend on the particular choice of $J$. \vv \qed
%
%
%
%
\pd \dd 2.11*. $ATh^n(\be_i,\ (C,\bb A))$ for\ $\be_i\in J^*$ is:
$Th^n(C,\bb A)|_{[\be_i,\ga)})$\ for every $\ga\in J^{**}$\ 
$\ga>\be_i$,\ $cf(\ga)=\om$; \ Where\ $J^{**}$\ is from Lemma 2.10*.  
(Actually this is $s_i$ from notation 2.17). \vv
\pd Remark. Again, fixing $J^*$ and $h$ it is easily seen that the definition 
does not depend on the choice of $J^{**}$.
%
%
%
%
\pd \dd 2.13*. Let \ $cf(C)>\om$, \ $M = (C,\bb A)$\  and we define
the model \ $g^n(M) = \bigr(C,g^n(\bb A)\bigl)$. \par
\\Let:\ \ $\bigl(g^n(\bb A)\bigr)_s = \{ \al_i<\la : 
s=ATh^n(\be_i,(C,\bb A)), \be_i\in J^*, h(\be_i)=\al_i \}$   (so this is   
a subset of $\la$). \par
\\let $d=l(\bb A)$\ and \ $T(n,d)$\ := the set of formally possible \ 
$Th^n(M,\bb B)$, where $l(\bb B)=d$. \par
We define a finite sequence of subsets of $\la$:
$$g^n(\bb A) := \lan\ldots, \bigr(g^n(\bb A)\bigl)_s,\ldots
\ran_{s\in T(n,d)}$$  \vv
%
%
%
%
\pt Lemma 2.14*. The analogs of lemma 2.14 hold for $g^n(C,\bb A)$  \vv
%
%
%
%
\pt Theorem 2.15*. If \ $cf(C)>\om$,\  then for each $n$ there
is an $m=m(n)$ such that if:  \par
\\$t_0=Th^m(C,\bb A)|_{\be_0}$, \ $t_1=WTh^m\bigl(C,g^m(C,\bb A)\bigr)$,
\ $t_2=ATh^m\bigl(\be_0,(C,\bb A)\bigr)$  \par
then we can effectively compute \ $Th^n(C,\bb A)$\ from\ $t_0,t_1,t_2$. 
(If \ $C$ has a first element $\de$, set $\be_0=\de$ and we don't need $t_0$).
\vv
\pd Remark. Following our notations, $Th^n(C,\bb A)$ is equal to 
$t_0+\sum_{i<\la}s_i$. By 2.10* we get for example 
(if $J^*=J^{**}$ from 2.10): 
$\sum_{i\le k<j}s_i = s_i$ for $cf(j)\le\om$ . \par
\\What we say in 2.13* is that if we know $t_0$ and $s_0$ and we know, 
roughly speaking, `how many' theories of every kind appear in the sum (this 
information is given by $t_1$), then we can compute the sum of the theories
exactly as in the case of well ordered chains. \vv
%
%
%
%
\pd Notation 2.16*. We will denote \ $ \langle t_0,t_1,t_2 \rangle $\ 
from Thm. 2.15* by \ $WA^{m(n)}$. \vv
\sec 
%
%
%
%
\hh 3. Major segments \par
In this section we define a theory \ $T$ \ which is interpretable in Peano 
arithmetic and reduce a supposed interpretation of \ $T$ \ in a chain $C$ to 
an interpretation of even a simpler theory in a chain $D$ having some  
favorable properties which will lead us to a contradiction.  \par
%
%
%
%
\pd \dd 3.0. Let $T$ be a first order theory with a signature  
consisting of one binary predicate $p$. The axioms of $T$ are as follows:

(a) $\all x\ex y\all z[p(z,y)\leftrightarrow z=x]$ \par
(b) $\all x\all y\ex u\all z[p(z,u)\leftrightarrow(p(z,x)\vee p(z,y))]$ \par
(c) $\ex x\all y[\neg p(y,x)]$ \vv
\\Intuitively (a) means that for every set $x$ there exists the set $\{x\}$,
\ (b) means that for every set $x,y$ there exists the set $x\cup y$ \ and
(c) means that the empty set (or an atom) exists. \par
\\Now, Peano arithmetic easily interprets $T$ in the sense of definition 1.4  
(let $d=1$, $U(x)$:=$x=x$, $E(x,y)$:=$x=y$ and $p'(x,y)$:= ``there exists a 
prime number $p$ \st $p^x$ divides $y$ but $p^{x+1}$ does not''),
so it suffices to show that no chain $C$ interprets $T$. \par
So Suppose $C$ is a chain that interprets $T$ by: 
$$I = \langle d,U(\bb X_1,\bb W),E(\bb X_1,\bb X_2,\bb W),
P(\bb X_1,\bb X_2,\bb W) \rangle .$$
\\We may assume, by changing $E$, that the interpretation is universal, 
i.e. $C\models(\forall\bb X)U(\bb X)$, and that the relation $P$ 
satisfies extensionality. $\bb W\sb C$ is a finite sequence of parameters 
and we will usually forget to write them.  Remember, for later stages, that  
we may assume that there is a closed cofinal subchain in $C$, if not add the 
completion of some cofinal subchain to $C$ and to the parameters and,  
if necessary, modify $I$. \par
\\Hence, the interpretation defines a model of $T$:
$${\cal M} = \lan \ ({\cal P}(C)^d/E), P \ \ran$$
\pd Notation. We will refer to ($d$-tuples of) subsets of $C$ as 
`elements'. If not otherwise mentioned, all the sequences appearing in the
formulas have length $d$ (= the dimension of the interpretation).   \par
\\We write $\bb X\sim\bb Y$ when ${\cal M}\models E(\bb X,\bb Y)$  \par
\\We write, for example, $\bb A\cap\bb B$ meaning $\langle 
A_0\cap B_0,\ldots,A_{lg(\bb A)-1}\cap B_{lg(\bb B)-1} \rangle$ and assuming 
$lg(\bb A)=lg(\bb B),\ \bb A = \langle A_0\ldots,A_{lg(\bb A)-1} \rangle$ etc.
\par
\\We also write $\bb A\subseteq C$ when $\bb A \in {\cal P} (C)^{lg \bb A}$. \par
%
%
%
%
\pd \dd 3.1. \par
1) A subchain $D\subseteq C$ is a {\sl segment} if it is convex (i.e. 
$x<y<z \ \&\ x{,}z\in D\ \Rightarrow \ y\in D$). \par
2) We will write $\bb A\sim\bb B$ when $\bb A,\bb B\subseteq C$ and 
$C'\models E(\bb A,\bb B)$.                          \par
3) Let $\bb A,\bb B\subseteq C$. We will say that $\bb A,\bb B$ 
{\sl coincide on ({\rm resp.} outside)} a segment $D\subseteq C$, 
if $\bb A\cap D=\bb B\cap D$ (resp. $\bb A\cap(C-D)=\bb B\cap(C-D) \ )$.  \par
4) The {\sl bouqet size} of a segment $D\sb C$ is the supremum 
of cardinals $|S|$ where $S$ ranges over collections of nonequivalent 
elements coinciding outside $D$.                             \par
5) A {\sl Dedekind cut} of $C$ is a pair $(L,R)$ where $L$ is an initial 
segment of $C,$\ $R$ is a final segment of $C$ and $L\cap R = \emptyset, 
L\cup R = C$.                                                    \vv
\\Our next step is to show that the bouquet size of every initial segment 
is either infinite or a-priory bounded.
%
%
%
%
\pt Lemma 3.2. There are monadic formulas $\theta_1(\bb X,\bb Z)$, and
$\theta_2(\bb X,\bb Y,\bb Z)$  such that: \par
1) For every finite, nonempty collection $S$ of elements, there is 
an element $\bb W$ such that for an arbitrary element $\bb A$, \ $C \models 
\theta_1(\bb A,\bb W)\ $ if and only if there is an element $\bb B\in S$ 
such that $\bb B\sim\bb A$.                   \par
2) For every finite, nonempty collection $S$ of pairs of elements, 
there is an element $\bb W$ such that for an arbitrary pair of elements 
$\langle\bb A_1,\bb A_2\rangle$,\ \ $C \models 
\theta_2(\bb A_1,\bb A_2,\bb W) \ $ if and only if there is a 
pair $\langle\bb B_1,\bb B_2\rangle \in S$ such that $\bb B_1\sim\bb A_1,
\bb B_2\sim\bb A_2$.                              \vv
\pd proof. Easy ($T$ allows coding of finite sets and $C$ interprets $T$).
\vv \qed
\\Thinking of $P$ as the $\epsilon$ relation, we will sometimes denote by 
something like $\big\{ \{\bb X\},\{\bb Y\},\ldots \big\}$ 
the set that codes $\bb X,\bb Y,\ldots$.
%
%
%
%
\pt Proposition 3.3. Fix a large enough $m<\omega$, (e.g. such that from 
$Th^m(C,\bb A_1,\bb A_2,\bb W)$ we can compute whether $C\models  
\theta_2(\bb A_1,\bb A_2,\bb W)$).\  \  Let: 
$N_1 = \big| \{ Th^m(D,\bb X,\bb Y,\bb Z) : D$\ is a chain, 
$\bb X,\bb Y,\bb Z \subseteq D\} \big |$   \par
{\sl Then,} for every Dedekind cut $(L,R)$ of $C$, either the bouquet size of
$L$ is at most $N_1$ and the bouquet size of $R$ is infinite, or, the
bouquet size of $R$ is at most $N_1$ and the bouquet size of $L$ is infinite.
\vv
\pd Proof. See \gu \ Thm. 6.1 and Lemma 8.1.
\vv \qed
%
%
%
%
\pd \dd 3.4. \par
1) A segment $D\subseteq C$ is called {\sl minor} if it's bouquet size is at 
most $N_1$.      \par
2) A segment $D\subseteq C$ is called {\sl major} if it's bouquet size is  
infinite.            \vv
%
%
%
%
\pt Conclusion 3.5. $C$ is major and for every Dedekind cut $(L,R)$ of $C$, 
either $L$ is minor and $R$ is major, or vice versa. \vv
\pd Proof. By Prop. 3.3. (and note that $T$  has only infinite models so
$C$ has an infinite number of $E$-equivalence classes).
\vv \qed
%
%
%
%
\pd \dd 3.6. An initial (final) segment $D$ is called a 
{\sl minimal major segment} if $D$ is major and for every initial (final) 
segment $D'\subset D$, \ $D'$ is minor.    \vv
%
%
%
%
\pt Lemma 3.7. There is a chain $C^*$ that interprets $T$ and an initial 
segment $D\subseteq C^*$ (possibly $D=C^*$) such that $D$ is a minimal major 
segment.   \vv
\pd Proof. (By \gu\ lemma 8.2). Let $L$ be the union of all the minor initial   
segments (note that if $L$ is minor and $L'\subseteq L$ then $L'$ is minor).
If L is major then set $L = D$ and we are done. Otherwise, let $D$ = $C-L$,
and by conclusion 3.5 $D$ is major. If there is a final segment 
$D'\subseteq D$ which is major then $C-D'$ is minor. But, $C-D'\supset L$, 
a contradiction.  \par
\\So $D$ is a minimal major (final) segment. Now take $C^{{\rm INV}}$ to be 
the inverse chain of $C$. By virtue of symmetry $C^{{\rm INV}}$ interprets $T$ 
and $D$ is a minimal major initial segment of $C^{{\rm INV}}$. 
\vv \qed
\pd Notation. Let $D\subseteq C$ be the minimal major initial segment we 
found in the previous lemma.   \vv
\pd Discussion. It is clear that $D$ is definable in $C$. (It's the shortest
initial segment such that there at most $N_1$ nonequivalent elements 
coinciding outside it).
What about $cf(D)$? It's easy to see that $D$ does not have a last  
point. On the other hand, it was proven in \gu\ that $T$ is not interpretable 
in the monadic theory of 
short chains (where a chain $C$ is short if every well ordered subchain of
$C$ or $C^{{\rm INV}}$ is countable). But we don't need to assume that 
the interpreting chain is short in order to apply \gu's argument. All we have
to assume, to get a contradiction is that $cf(D)=\omega$ (which is of course  
the only possible case when $C$ is short).
So, if $C$ interprets $T$ and $cf(D)=\omega$, we can repeat the argument 
from \gu\ to get a contradiction. Therefore, we can conclude:  \vv
%
%
%
%
\pt Proposition 3.9. $cf(D) > \omega$   \vv \qed
%
%
%
%
\pd Notation 3.10. \tk\ will denote the theory of a family of $k$ sets and 
the codings of every subfamily. \vv
\pd Discussion (continued).  Now, fix an element $\bb R\subseteq (C-D)$ 
witnessing the fact that $D$ is major, and define:
$$ S = \{\ \bb A\subseteq C :\  \bb A\cap (C-D) = \bb R\ \} $$
\\So $S/_\sim$ is infinite by the choice of $\bb R$ (and of course definable 
in $C$ with an additional parameter $\bb R$). For the moment let $k=2$ and
fix a finite subset of 6 nonequivalent elements in $S$, 
$\lan \bb A_1,\bb A_1\ldots \bb A_6 \ran$. \ We want to define in $D$ a 
structure that contains 2 `atoms' and 4 codings by using the $\bb A_i$'s. \par
\\Since ${\cal M}\models T$ \ we have an element $\bb W\sb C$ 
(not necessarily in $S$) which can be identified with the set: \par
$\big\{ \{\bb A_1\},\{\bb A_2\},\{\{\bb A_3\}\},\{\{\bb A_4,\bb A_1\}\}, 
\{\{\bb A_5,\bb A_2\}\},\{\{\bb A_6,\bb A_1\} \{\bb A_6,\bb A_2\}\} \big\}$.

\\Look at the following formulas: \par
Atom$(\bb X,\bb W)$ := $P(\{\bb X\},\bb W)$  \par
Set$(\bb Y,\bb W)$ := $\neg{\rm Atom}(\bb Y,\bb W)$ \&  
$\ex\bb Z\ex\bb V(P(\bb Z,\bb W)\&P(\bb V,\bb Z)\&(P(\bb Y,\bb V))$  \par
Code$(\bb X,\bb Y,\bb W)$ := Atom$(\bb X,\bb W)$ \& Set$(\bb Y,\bb W)$ \& 
$(\ex\bb Z(P(\bb Z,\bb W)\&P(\{\bb X,\bb Y\},\bb Z))$  \par
\\Using these formulas we can easily define in $C$ a structure which 
satisfies $T_2$, where $\bb A_1$ and $\bb A_2$ are the atoms $\bb A_3$ 
codes the empty subfamily, $\bb A_4$ codes $\bb A_1$ etc. \  
But for every natural number \thinspace $k$ \thinspace we can define a 
structure for \tk\ by picking $k+2^k$ elements from $S$ and a suitable 
$\bb W$, and note that the above formulas do not depend on k.  \par
Now we claim that we can interpret \tk\ even in $D$ and not in all $C$. To see  
that, look at the formula \ Code$(\bb X,\bb Y,\bb W)$. 
There is an $n<\om$ \st we can decide from $Th^n(C,\bb X,\bb Y,\bb W,\bb R)$ 
if \ Code$(\bb X,\bb Y,\bb W)$ holds. By the composition theorem it suffices 
to look at 
$Th^n(D,\bb X\cap D,\bb Y\cap D,\bb W\cap D,\bb R\cap D)$ \ and \ 
$Th^n(C-D, \bb X\cap (C-D),\bb Y\cap (C-D),\bb W\cap (C-D),\bb R\cap (C-D))$. 
But, since we restrict ourselves only to elements in $S$, the second theory 
is constant for every $\bb X,\bb Y$ in $S$. It is: \ 
$Th^n(C-D, \bb R,\bb R,\bb W\cap (C-D),\bb R)$. So it suffices to know only 
$Th^n(D,\bb X\cap D,\bb Y\cap D,\bb W\cap D)$, ($\bb R\cap D = \emptyset$). 
Now use Lemma 2.4 to get a formula Code*$(\bb X,\bb Y,\bb W\cap D)$ that 
implies Code$(\bb X\cup\bb R,\bb Y\cup\bb R,\bb W)$, and the same holds for 
the other formulas (including the equality formula for members of $S$). \par
We get an interpretation of \tk\ on $D$ with an additional parameter $\bb W$. 
Remember that we allowed parameters $\bb V$ in the original interpretation 
of $T$ in $C$ and we can assume that $\bb W$ is a sequence that contains the  
coding set  and the old parameters (all intersected with $D$). \par
\\The universe formula of the interpretation is 
Atom*$(\bb X,\bb W)$ $\vee$ Set*$(\bb X,\bb Y,\bb W)$, the coding formula 
is Code*$(\bb X,\bb Y,\bb W)$ and the equality formula is 
$E^*(\bb X,\bb Y,\bb W)$.  And for different $k$'s and even different choices 
of members of $S$, the formulas (and their quantifier depth) are unchanged
except for the parameters $\bb W$. \par
\\It is easy to see that, since $D$ is minimal major, for every proper 
initial segment $D'\sbb D$ there are no more then $N_1$ (from definition 3.4) 
$E^*$ nonequivalent members of $S$ coinciding outside $D'$. We will say,
by abuse of definition, that $D$ is still a minimal major initial segment 
with respect to $E^*$. \ To sum up, we have proven:  
%
%
%
%
\pt Theorem 3.11. If there is an interpretation of $T$ in the monadic theory
of a chain $C$ \ then, there is a chain $D$ such that $cf(D)>\om$, and 
\st for every $k<\om$ there is an interpretation of \tk\ in the monadic 
theory of $D$ such that the interpretation does not 
``concentrate'' on any proper initial segment of $D$ 
(i.e. $D$ itself is the minimal major initial segment of $D$).
Furthermore, there is an $n<\om$  which does not depend on k, \st all the 
interpreting formulas have quantifier depth $<n$.
\vv \qed
\sec
%
%
%
%
%
\hh 4.Preservation of theories under shufflings \par  
We will 
define here shufflings of subchains and show that the partial theories defined 
in \S 2 are preserved under them. \medbreak
\\{\bf Convention:} 1. Throughout this section, $\de$ will denote an 
ordinal with $cf(\de)=\la>\omega$. \par
\\2. Unless otherwise said, all the chains mentioned in this section
are well ordered chains (i.e. ordinals).  We will deal with general chains in 
the next section. \medbreak
%
%
%
%
\pd \dd 4.1.  \par
\\1) Let $a\sb\la$. We say that $a$ is a {\sl \sc\ subset of} $\la$ if for 
every $\al<\la$ with $cf(\al)>\om$:  \par
\\if $\al\in a$ then there is a club subset of $\al,\ C_\al$ \st $C_\al\sb a$ 
\ and if $\al\not\in a$ then there is a club subset of 
$\al,\ C_\al$ \st $C_\al\cap a=\emptyset$.  \par
\\Note that $\la$ and $\emptyset$ are semi-clubs
and that a club $J\sb\la$ is a \sc\ provided that the first and the
successor points of $J$ are of cofinality $\leq\omega$. \par
\\2) Let $X,Y\sb\de, \ J =\{\al_i:i<\la\}$ a club subset of $\de$, and let
$a\sb\la$ be a \sc\ of $\la$. We will define the {\sl shuffling of
$X$ and $Y$ with respect to $a$ and $J$},\ denoted by $[X,Y]^{J}_a$, as:
  $$[X,Y]^{J}_a = \bigcup_{i\in a}\big(X\cap[\al_i,\al_{i+1})\big)\cup
                     \bigcup_{i\not\in a}\big(Y\cap[\al_i,\al_{i+1})\big)$$
\\3) When $J$ is fixed (which is usually the case), we will denote the 
shuffling of $X$ and $Y$ with respect to $a$ and $J$,\ by $[X,Y]_a$.  \par
\\4) When $\bb X,\bb Y\sb\de$ are of the same length, we define 
$[\bb X,\bb Y]_a$ naturally.                                    \par
\\5) We can define shufflings naturally when $J\sbb\de$ is a club, 
and $a\sb otp(J)$ is a \sc.     \vv
%
%
%
%
\pd Notation 4.2.    \par
\\1) Let $\bb P_0\sb\de$ and $J\sb\de$ a club subset of $\de$ witnessing 
$ATh(\de,\bb P_0)$ as in lemma 2.10. For $n<\om$, and $\be<\ga$ with 
$\ga\in J, cf(\ga)=\om$, 
we denote $Th^n\big((\de,\bb P_0)|_{[\be,\ga)\big)} = 
ATh^n\big(\be,(\de,\bb P_0)\big)$ by $s^{n}_{\bb P_0}(\be)$ or just 
$s^{n}_{0}(\be)$.
(Of course, this does not depend on the choice of $J$ and $\ga$).  \par
\\2) When $n$ is fixed we denote this theory by $s_{\bb P_0}(\be)$ or 
$s_0(\be)$.                                                        \par
\\3) Remember: $g^n(\bb P_0)_s$ is the set $\{\be<\de : 
s^{n}_{\bb P_0}(\be)=s\}$. (See def. 2.13.)    \par
\\4) $S_{0}^{\de}$ is the set $\{\ga<\de : cf(\ga)=\om\}$.   \vv
%
%
%
%
\pd \dd 4.3. Let $\bb P_0,\bb P_1\sb\de$ be of the same length and $J\sb\de$
be a club. We will say that $J$ is $n$-{\sl suitable for} $\bb P_0,\bb P_1$
if the following hold:   \par
a) $J$ witnesses $ATh(\de,\bb P_l)$ for $l=0,1$.           \par
b) $J=\{\al_i:i<\la\}, \al_0=0$ and $cf(\al_{i+1})=\om$.       \par
c) $J\cap g^n(\bb P_l)_s\cap S_{0}^{\de}$ is either a stationary subset of 
$\de$ or is empty.                                                 \par
\\When $n\geq 1$ and $WA^n(\de,\bb P_0)=WA^n(\de,\bb P_1)$             
(see notation 2.15) we require also that:                              \par
d) If $\al_j\in J$ $cf(\al_j)\le\om$ and $s_l(\al_j)=s$ then there are 
$k_1,k_2<\om$ \st\ $s_l(\al_{j+k_1})=s,$ and $s_{1-l}(\al_{j+k_2})=s$. \vv 
\pd Remark. It is easy to see that for every finite sequence
 
$\lan \bb P_0,\bb P_1,\ldots,\bb P_n \ran \sb\de$ with equal lengths, there  
is a club $J\sb\de$ which is $n$-suitable for every pair of the $\bb P_i$'s.
\medbreak
\\We will show now that $ATh$ is preserved under `suitable' shufflings.
%
%
%
%
\pt Theorem 4.5. Suppose that $\bb P_0,\bb P_1\sb\de$ are of the same length,
$n\geq 1$ and $WA^n(\de,\bb P_0)=WA^n(\de,\bb P_1)$. 
(In particular, $ATh^n(0,(\de,\bb P_0))=ATh^n(0,(\de,\bb P_1)) := t $).
Let $J\sbb\de$ be $n$-suitable for $\bb P_0,\bb P_1$ of order type $\la$ and
$a\sb \la$ a \sc.
Then, $ATh^n(0,(\de,[\bb P_0,\bb P_1]^{J}_{a}))=t$   \vv
\pd Proof. Denote $\bb X:=[\bb P_0,\bb P_1]^{J}_{a}$.  \par
\\We will prove the following facts by induction on $0<j<\la$:       \par
\\$(*)$  For every $i<j<\la$ with $cf(j)\leq\om$:                  \par
$i\in a \imp Th^n([\al_i,\al_j),\bb X) = 
Th^n([\al_i,\al_j),\bb P_0) = s_0(\al_i).$                          \par
$i\not\in a \imp Th^n([\al_i,\al_j),\bb X) = 
Th^n([\al_i,\al_j),\bb P_1) = s_1(\al_i).$                          \par
\\$(**)$  For every $i<j<\la$ with $cf(j)>\om$:                     \par
$i,j\in a \imp Th^n([\al_i,\al_j),\bb X) = 
Th^n([\al_i,\al_j),\bb P_0).$                                       \par

$i,j\not\in a \imp Th^n([\al_i,\al_j),\bb X) = 
Th^n([\al_i,\al_j),\bb P_1).$                                        \par\par
\\In particular, by choosing $i=0$ we get (remember $\al_0=0$),
$Th^n([0,\al_j),\bb X) = t$ whenever $cf(\al_j)=\om$.                \par\par
\\\underbar{$j=1$} (so $i=0$): \ Let $l=0$ if $i\in a$ and $l=1$ if 
$i\not\in a$. So $\bb X\cap[0,\al_j)=\bb P_l\cap[0,\al_j)$ and so 
$Th^n([0,\al_j),\bb X) = Th^n([0,\al_j),\bb P_l) = t$                 \par\par
\\\underbar{$j=k+1<\om$}: There are 4 cases. Let us check for example the 
case $i\in a, \ j-1=k\not\in a$. \ By the composition theorem (2.7) and the 
induction hypothesis we have:                                          \par
\\$Th^n([\al_i,\al_j),\bb X)$ \  = \  
$Th^n([\al_i,\al_k),\bb X)+Th^n([\al_k,\al_{k+1}),\bb X)$ \  = \ 
$s_0(\al_i)+Th^n([\al_k,\al_j),\bb P_1)$ \ = \ $s_0(\al_i)+s_1(\al_k).$
So we have to prove $s_0(\al_i)+s_1(\al_k) = s_0(\al_i)$.               \par
\\Since $J$ is $n$-suitable there is an $m<\om$ \st 
$s_0(\al_{i+m})=s_1(\al_k)$ and so,                                     \par
\\$ s_0(\al_i) \ = \ Th^n([\al_i,\al_{i+m+1}),\bb P_0) \ = \ 
Th^n([\al_i,\al_{i+m}),\bb P_0)+Th^n([\al_{i+m},\al_{i+m+1}),\bb P_0) \  = \ 
s_0(\al_i)+s_0(\al_{i+m}) \ = \ s_0(\al_i)+s_1(\al_k).$                 \par
\\So $s_0(\al_i)+s_1(\al_k) = s_0(\al_i)$ as required.                  \par
The other cases are proven similarilly.                              \par\par
\\\underbar{$j=\om$}: Suppose $i<\om,i\in a$. We have to prove that
$Th^n([\al_i,\al_\om),\bb X) = s_0(i)$. Now either $(\la\setminus a)\cap\om$  
is unbounded or $a\cap\om$ is unbounded and suppose the first case holds. 
Let $i<i_0<i_1\ldots $ be a strictly increasing sequence in $(\la\setminus a)\cap\om$ .
By the induction hypothesis we have:                                \par
\\$Th^n([\al_i,\al_\om),\bb X) \ = \ Th^n([\al_i,\al_{i_1}),\bb X)+\sum_
{0<m<\om} Th^n([\al_{i_m},\al_{i_{m+1}}),\bb X) \ = \ s_0(\al_i)+\sum_
{0<m<\om}s_1(\al_{i_m}).$                                          \par
Now choose (using the suitability of $J$), a strictly increasing sequence
$\be_{i_0}<\be_{i_1}\ldots\sb\la$ \st $\be_{i_m}=\al_{j_m+1}$ for some 
$j_m<\la$, $\be_{i_1}>\al_i$ and 
\st for every $0<m<\om, \ s_0(\be_{i_m})=s_1(\al_{i_m})$.
We will get:                                                           \par
\\$s_0(\al_i) \ = \ Th^n([\al_i,\al_\om),\bb P_0)$ 
\ = \  $Th^n([\al_i,\be_{i_1}),\bb P_0)+\sum_{0<m<\om}
                        Th^n([\be_{i_m},\be_{i_{m+1}}),\bb P_0)$  \ =   \par
\\$s_0(\al_i)+\sum_{0<m<\om}s_0(\be_{i_m})$ 
\ = \  $s_0(\al_i)+\sum_{0<m<\om}s_1(\al_{i_m})$.                     \par
\\So we have $s_0(\al_i) = Th^n([\al_i,\al_\om),\bb X)$ as required.  \par
\\When only the other case holds (i.e. only $a\cap\om$ is unbounded) the 
proof is easier.                                                      \par
\\When $i\not\in a$ we prove similarly  that $Th^n([\al_i,\al_\om),\bb X) = 
s_1(\al_i)$                \par\par
\\\underbar{$cf(j)=\om$}: \ Choose a sequence (in $a$ or $\la\setminus a$), 
$i<i_0<i_1\ldots\ \ \sup_{m}i_m=j, \ i_m$ non limit, and continue as in the 
case  $j=\om$.                     \par\par
\\\underbar{$cf(j)>\om$}: \ Now we have to check $(**)$. \par
\\So suppose $i,j\in a$ and we have to show $Th^n([\al_i,\al_j),\bb X) = 
Th^n([\al_i,\al_j),\bb P_0)$.                                \par
\\Let $\{\be_\ga : \ga<cf(j)\}\sb a$ be a club subset of $j$ with 
$\be_0=i$. By the induction hypothesis we have: 
$Th^n([\al_i,\al_j),\bb X) =
\sum_{\ga<cf(j)}Th^n([\be_\ga,\be_{\ga+1}),\bb X) =
\sum_{\ga<cf(j)}s_0(\be_\ga) =
\sum_{\ga<cf(j)}Th^n([\be_\ga,\be_{\ga+1}),\bb P_0)$ = 
$Th^n([\al_i,\al_{j}),\bb P_0)$ as required.                  \par
\\The case $i,j\not\in a$ is similar.                          \par
\\\underbar{$j=k+2$}: \ Easy.                                 \par\par
\\\underbar{$j=k+1, cf(k)=\om$}: \ Easy.                      \par\par
\\\underbar{$j=k+1,cf(k)>\om$}: There are 8 cases. We will check for example
the case: $0<i$, $i\in a$, $k\not\in a$.                              \par
\\Choose $\{i_\ga : \ga<cf(k)\}\sb\la\setminus a$ \ a club \st $i<i_0$ and 
$s_0(\al_i)=s_1(\al_{i_0})$.                                 \par
\\Note that ($i\not=0$) 
$s_0(\al_i)+s_0(\al_i)=s_0(\al_i)=s_0(\al_i)+s_1(\al_{i_0})$.   \par
\\So we get \ $Th^n([\al_i,\al_j),\bb X)$ = 

\\$Th^n([\al_i,\al_{i_0}),\bb X)+
\sum_{\ga<cf(k)}Th^n([\al_{i_\ga},\al_{i_{\ga+1}}),\bb X)+
Th^n([\al_k,\al_j),\bb X)$ =                                  \par
\\$s_0(\al_i)+\sum_{\ga<cf(k)}s_1(\al_{i_\ga})+s_1(\al_k) = 
s_0(\al_i)+s_1(\al_{i_0})+\sum_{0<\ga<cf(k)}s_1(\al_{i_\ga})+s_1(\al_k) = 
s_0(\al_i)+s_0(\al_i)+\sum_{0<\ga<cf(k)}s_1(\al_{i_\ga})+s_1(\al_k) = 
s_0(\al_i)+\sum_{0<\ga<cf(k)}s_1(\al_{i_\ga})+s_1(\al_k)$ =    \par
\\$s_1(\al_{i_0})+\sum_{0<\ga<cf(k)}s_1(\al_{i_\ga})+s_1(\al_k) = 
\sum_{\ga<cf(k)}s_1(\al_{i_\ga})+s_1(\al_k) = 
Th^n([\al_i,\al_{i_0}),\bb P_1) =
Th^n([\al_i,\al_j),\bb P_0) = s_0(\al_i)$.
So $Th^n([\al_i,\al_j),\bb X) = s_0(\al_i)$ as required.           \par
\\Check the other cases: when $i=0$ use the fact $s_0(0)=s_1(0)$.  \par\par
So we have gone through all the cases and proven $(*)$ and $(**)$.  \vv \qed
%
%
%
%
\pt Conclusion 4.6. Let $\bb P_0,\bb P_1\sb\de$, $WA^n(\de,\bb P_0) = 
WA^n(\de,\bb P_1)$, $a\sb\la$ a \sc\ and $J\sb\de$ an $n$-suitable club subset
for $\bb P_0$,$\bb P_1$. Then:    \par
1) $J$ is an $n$-suitable club subset for the pair  
$\bb P_0$, $[\bb P_0,\bb P_1]^J_a$.       \par
2) $[g^n(\bb P_0),g^n(\bb P_1)]^J_a\cap J = g^n([\bb P_0,\bb P_1]^J_a)\cap J$
\ (so $[g^n(\bb P_0),g^n(\bb P_1)]^J_a$ and $g^n([\bb P_0,\bb P_1]^J_a)$
have the same $WTh$). \vv
\pd Proof. Use $(*),(**)$ from the last theorem.    \vv \qed
Our next aim is to show that $WTh$ (hence, by 2.15 and 4.6(2), also $Th$) 
is preserved under shufflings.
%
%
%
%
\pd \dd 4.7. Let $a,\bb P\sb\la$. We define $a-WTh^n(\la,\bb P)$ by induction
on $n$:           \par
for $n=0$: \ \ $a-WTh^0(\la,\bb P) = \{ t : th(\la,(\bb P,a)) {\rm\ is\  
stationary\ in\ } \la \}$  \ \ (see def. 2.9)   \par
for $n+1$: \ \ $a-WTh^{n+1}(\la,\bb P) = \{ \lan S_1^{\bb P,a}(Q), 
S_2^{\bb P,a}(Q),S_3^{\bb P,a}(Q) \ran : Q\sb\la \}$  \ \ Where:  \par
\\$S_1^{\bb P,a}(Q) = a-WTh^n(\la,\bb P,Q)$                \par
\\$S_2^{\bb P,a}(Q) = \big\{ \lan t,s\ran : 
\{\be\in a: WTh^n(\la,\bb P,Q)|_\be=t,\ th(\be,\bb P,Q)=s\}\ {\rm is\   
stationary\ in\ }\la \big\}$                     \par
\\$S_3^{\bb P,a}(Q) = \big\{ \lan t,s\ran : 
\{\be\in \la\setminus a: WTh^n(\la,\bb P,Q)|_\be=t,\ th(\be,\bb P,Q)=s\}\ {\rm is\ 
stationary\ in\ }\la \big\}$                         \vv
%
%
%
%
\pd Remark 4.8. 0) Remember that if $\bb P\sb\de$ and $J\sb\de$ is a club, 
then  $WTh^n(\de,\bb P\cap J)=WTh^n(\de,\bb P)$. 
Moreover, if $J\sb\de$ club of order type $\la$ and 
$h\colon J \rightarrow \la$ is the isomorphism between $J$ and $\la$, 
then for every $\bb P\sb\de$, $WTh^n(\de,\bb P) = WTh^n(\la,h(\bb P\cap J))$.

\\1) $WTh^n(\la,\bb P)$ tells us if certain sets are stationary. 
$a-WTh^n(\la,\bb P)$ tells us if their intersection with $\la$ and $\la\setminus a$ 
are stationary.                                           \par
\\2) We could have defined $a-WTh^n(\la,\bb P)$ by $WTh^n(\la,\bb P,a)$,
which gives us the same information. We prefared the original definition
because it seems to be easier to see the preservation under shufflings using 
it. \vv
%
%
%
%
\pt Fact 4.9. For any $a\sb\la$, $WTh^n(\la,\bb P)$ is effectively computable 
from $a-WTh^n(\la,\bb P)$, so if $\bb P,\bb Q\sb\la$ and 
$a-WTh^n(\la,\bb P)=a-WTh^n(\la,\bb Q)$ 
then:  $WTh^n(\la,\bb P)=WTh^n(\la,\bb Q)$.  \vv
\pd Proof. Trivial.                          \vv \qed
%
%
%
%
\pt Theorem 4.10. Suppose $a,J,\bb P_0,\bb P_1\sb\la, \ a$ \sc, $J$ club,  
$\bb X:= [\bb P_0,\bb P_1]^J_a$ and  
$a-WTh^n(\la,\bb P_0)=a-WTh^n(\la,\bb P_1)$.    \par
\\Then: $a-WTh^n(\la,\bb P_0)=a-WTh^n(\la,\bb X)$. (It follows
$WTh^n(\la,\bb X)=WTh^n(\la,\bb P_0)=WTh^n(\la,\bb P_1)$ ).  \vv
\pd Proof. by induction on $n$  (for every $a',J',\bb X',\bb Y'$):   \par
\\$n=0$: \ Check.                                                    \par
\\$n+1$: \ Suppose $Q_0\sb\la$ and $\{ \lan S_1^{\bb P_0,a}(Q_0),    
S_2^{\bb P_0,a}(Q_0),S_3^{\bb P_0,a}(Q_0) \ran \} \in 
a-WTh^{n+1}(\la,\bb P_0)$. 
Choose (using the equality of the theories) $Q_1\sb\la$ \ \st 

\\$\{ \lan S_1^{\bb P_1,a}(Q_1), S_2^{\bb P_1,a}(Q_1),S_3^{\bb P_1,a}(Q_1) 
\ran \} \in a-WTh^{n+1}(\la,\bb P_1)$, and \st the two triples are equal.
Define $Q_X := [\bb Q_0,\bb Q_1]^J_a$.                               \par
\\By the \ind $a-WTh^n(\la,\bb P_0,Q_0) = a-WTh^n(\la,\bb X,Q_X)$ so 
$S_1^{\bb P_0,a}(Q_0) = S_1^{\bb X,a}(Q_X)$.                          
Now suppose $\lan t,s \ran \in S_2^{\bb P_0,a}(Q_0), t\not=\emptyset$. 

\\Let $B^{\bb P_0}_{t,s} := \big\{ \be\in a : WTh^n(\la,\bb P_0,Q_0)|_\be=t,
\ th(\be,\bb P_0,Q_0)=s \big\}$ and this is a stationary subset of $\la$.
But for each such $\be$, since $t\not=\emptyset\imp cf(\be)>\om$, \ $a$
contains a club $C_\be\sb\be$ and, remembering a previous remark, we can  
restrict ourselves to $(\bb P_0,Q_0)\cap C_\be$.                        \par
\\Now suppose:                                                         
$a=\lan i_\ga : \ga<\la \ran$  (note that $a$ has to be stationary otherwise 
$S_2$ is empty) and $J=\lan \al_{i_\ga} : \ga<\la \ran$.                
Look at the club $J' = \lan \al_{i_\ga} : \al_{i_\ga}=i_\ga \ran$ and let
$J''$ := the accumulation points of $J'$.                              
\\Now $B^{\bb P_0}_{t,s}\cap J''$ is also stationary, and choose 
$\be$ in this set, and a club $C_\be\sb a\cap J'$.                    
By the choice of $C_\be$ we get:
$(\bb P_0,Q_0)\cap C_\be$ = $(\bb X,Q_X)\cap C_\be$, and this implies: 
$WTh^n(\la,\bb X,Q_X)|_\be=t$, and $th(\be,\bb X,Q_X)=s$. 
So, (since $\be$ was random) $B^{\bb X}_{t,s}$ is also stationary.   \par
\\The case $t=\emptyset$ is left to the reader.  We deal with $S_3$ 
symmetrically, replacing $a$ with $\la\setminus a$.     \par
\\So we have proven that $a-WTh^{n+1}(\la,\bb P_0)\sb a-WTh^{n+1}(\la,\bb X)$ \par
\\ Now, for the inverse inclusion suppose $Q_X\sb\la$ and:
$\lan S_1^{\bb X,a}(Q_X),S_2^{\bb X,a}(Q_X),S_3^{\bb X,a}(Q_X) \ran 
\in a-WTh^{n+1}(\la,\bb X)$.  Choose $R_0$ \st $R_0\cap a=Q_X\cap a$ and 
$R_1$ \st $R_1\cap(\la\setminus a)=Q_X\cap(\la\setminus a)$. 
\\Now choose $T_0$ \st $S_3^{\bb P_0,a}(T_0) = S_3^{\bb P_1,a}(R_1)$ and 
$T_1$ \st $S_2^{\bb P_0,a}(R_0) = S_2^{\bb P_1,a}(T_1)$.  
Let $Q_0$ be equal to $R_0$ on $a$ and to $T_0$ on $\la\setminus a$.
Let $Q_1$ be equal to $T_1$ on $a$ and to $R_1$ on $\la\setminus a$. 
It can be easily checked that: 
$\lan S_1^{\bb P_0,a}(Q_0),S_2^{\bb P_0,a}(Q_0),S_3^{\bb P_0,a}(Q_0) \ran$ = 
$\lan S_1^{\bb P_1,a}(Q_1),S_2^{\bb P_1,a}(Q_1),S_3^{\bb P_1,a}(Q_1) \ran$.
But $Q_X=[Q_0,Q_1]^J_a$, hence this triples are, by the same arguments as in
first part of the proof, equal to 
$\lan S_1^{\bb X,a}(Q_X),S_2^{\bb X,a}(Q_X),S_3^{\bb X,a}(Q_X) \ran$  \par
\\This proves the inverse inclusion: 
$a-WTh^{n+1}(\la,\bb P_0)\supseteq a-WTh^{n+1}(\la,\bb X)$, hence the equality  
\\$a-WTh^{n+1}(\la,\bb P_0)=a-WTh^{n+1}(\la,\bb P_1)=a-WTh^{n+1}(\la,\bb X)$
\vv\qed
%
%
%
%
\pd Notation 4.11. Suppose $\bb P,J\sb\de$, $J$ club of order type $\la$
and $a\sb\la$ a \sc.  \par
\\Let $t_1 := ATh^m(0,(\de,\bb P))$ and (keeping in mind 
remark 4.8.(0) ), let $h\colon J \rightarrow \la$ be the isomorphism 
between $J$ and $\la$ and let $t_2 := a-WTh^m(\la,h(g^m(\de,\bb P)\cap J))$ 
\par
\\We denote $\lan t_1,t_2 \ran$ by $a-WA^m(\de,\bb P)$  
(assuming $J$ is fixed).                                           \par
Collecting the last results we can conclude:
%
%
%
%
\pt Theorem 4.12. Let $J,\bb P_0,\bb P_1\sb\de$, $lg(\bb P_0)=lg(\bb P_1)$, 
$J$ an $n$-suitable club for $P_0,P_1$ of order type $\la$ and $a\sb\la$ a 
\sc\ and set $\bb X := [\bb P_0,\bb P_1]^J_a$.                     \par
\\Then: $a-WA^m(\de,\bb P_0) = a-WA^m(\de,\bb P_1) \imp 
a-WA^m(\de,\bb P_0) = a-WA^m(\de,\bb X)$, and in particular, if $m=m(n)$ then:    
$Th^n(\de,\bb P_0) = Th^n(\de,\bb P_1) = Th^n(\de,\bb X)$.        \vv
\pd Proof.  The first statement follows directly from 4.5 and 4.10. \par
\\For the second, by the definition of \enspace $a-WA$, and by 4.8(0), 4.9,  
equality of $a-WA^{m(n)}$ implies equality of $WA^{m(n)}$ from definition 2.16. 
But by 2.15 this implies the equality of $Th^n$.      \vv \qed
%
%
%
%
\sec
%
%
%
\hh 5. Formal shufflings   \par
In the previous section we showed how to shuffle subsets of well ordered 
chains and preserve their theories. Here we present the notion of formal 
shufflings in order to overcome two difficulties: \par
\\1. It could happen that the interpreting chain is of cofinality $\la$ but 
of a larger cardinality. Still, we want to shuffle objects of cardinality 
$\le\la$. The reason for that is that the contradiction we want to reach 
depends on shufflings of elements along a generic \sc\ added by the forcing, 
and a \sc\ of cardinality $\la$ will be generic only with respect to objects 
of cardinality $\le\la$. So we want to show now that we can shuffle theories, 
rather than subsets of our given chain.   \par
\\2. We want to generalize the previous results, which were proven for well 
ordered chains, to the case of a general chain. \medbreak
\pd Discussion. Suppose we are given a chain $C$ and a finite sequence of 
subsets $\bb A\sb C$ and we want to compute $Th^n(C,\bb A)$. As before we can 
choose an $n$-suitable club $J=\lan \al_i : i<\la \ran$ witnessing 
$ATh^n(C,\bb A)$ and letting $s_i:=Th^n(C,\bb A)|_{[\al_i,\al_{i+1})}$ we 
have: 
$Th^n(C,\bb A)=\sum_{i<\la}s_i$. \ Theorem 2.15 says that (for a large enough
$m=m(n)$ ) $WA^m(C,\bb A)$ which is $s_0$ and $WTh^m(\la,g^n(C,\bb A))$, 
determines $Th^n(C,\bb A)$. ( $g^n(C,\bb A)$ is a sequence of subsets of $\la$
of the form $g_s = \{ i : s_i=s\}$ ). \par
\\Moreover, since we have only finitely many possibilities for 
$WA^m(C,\bb A)$, we can decide whether $\sum_{i<\la}s_i = t$ inside
$H(\la^+) := \{ x : x$ is hereditarilly of cardinality smaller than $\la^+ \}$ 
even if the $s_i$'s are theories of objects of cardinality greater than 
$\la$. This motivates our next definitions: \par
%
%
%
%
\pd \dd 5.1. fix an $l<\om$
\item{1)} $S=\lan s_i : i<\la \ran$ is an {\sl \fps} if each $s_i$ is a 
formally possible member of 
$\{Th^n(D,\bb B) : D\ {\rm is\ a\ chain,\ }\bb B\sb D,\ lg(\bb B)=l \}$, and  
for every $i<j<\la$ with $cf(j)\le\om$ we have $s_i=\sum_{i\le k<j}s_k$.
\item{2)} The \fps\ $S$ is {\sl realized} in a model $N$ if there are
$J,C,\bb A$ as usual in $N$, and $s_i:=Th^n(C,\bb A)|_{[\al_i,\al_{i+1})}$.
\item{3)} Let $S=\lan s_i : i<\la \ran$, $T=\lan t_i : i<\la \ran$ be 
\fpss, $a\sb\la$ a \sc. We define the {\sl formal shuffling of $S$ and $T$
with respect to $a$ as: $[S,T]_a := \lan u_i : i<\la \ran$ where 
$$u_i=\cases{s_i &if $i\in a$ \cr t_i &if $i\not\in a$ \cr}$$  \vv
%
%
%
%
\pt Fact 5.2. 1. Let $\bb A,\bb B\sb C$ of length $l$, 
$J=\lan \al_i : i<\la \ran$ an $n$-suitable club and $a\sb\la$ a \sc.
Let $s_i:=Th^n(C,\bb A)|_{[\al_i,\al_{i+1})}$, $S=\lan s_i : i<\la \ran$, 
$t_i:=Th^n(C,\bb B)|_{[\al_i,\al_{i+1})}$, $T=\lan t_i : i<\la \ran$. \ \ 
Then: $S$ and $T$ are \fpss, and $[S,T]_a$ = 
$\lan Th^n(C,[\bb A,\bb B]_a^J)|_{[\al_i,\al_{i+1})} : i<\la \ran$.  \par
\\2. If in addition $WA^{m(n)}(C,\bb A)=WA^{m(n)}(C,\bb B)$, then $[S,T]_a$ is
an \fps. \par
\\3. If in addition $a-WA^{m(n)}(C,\bb A)=a-WA^{m(n)}(C,\bb B)$, then
$\sum_{i<\la}s_i=\sum_{i<\la}t_i=\sum_{i<\la}u_i$.  \vv
\pd Proof. Part 1 is obvious, part 2 follows from theorem 4.5  
and part 3 from 4.12. \vv \qed
%
%
%
%
We can define in a natural way the partial theories $WTh^m$ and $a-WTh^m$.\par
\pd \dd 5.3. For  $S=\lan s_i : i<\la \ran$ an \fps, denote $g^n(S)_s$:=
$\lan i<\la : s_i=s \ran$ and $g^n(S)$:=$\lan g^n(S)_s : s\ {\rm is\ a\ 
formally\ possible\ } n{\rm -theory} \ran$. \ \ We define $WTh^m(S)$ to be 
$WTh^m(\la,g^n(S))$, and for $a\sb\la$ a \sc, 
$a-WTh^m(S)$ is $a-WTh^m(\la,g^n(S))$. \par
\\Finally we define $a-WA^m(S)$ to be the pair $\lan s_0, a-WTh^m(S) \ran$. \vv
%
%
%
%
\pt Theorem 5.4. \ If \ $C,\bb A,J,S$\ are as usual then we can compute 
$Th^n(C,\bb A)$ from $WA^{m(n)}(S)$, moreover, the computation can be done in
$H(\la^+)$ even if $|C|>\la$. \vv
\pd Proof. The first claim is exactly 2.15. The second follows from the fact 
that $S$ and $WA^{m(n)}(S)$ are elements of $H(\la^+)$ and so is the 
correspondence between the (finite) set of formally possible $WA^m$'s and the 
formally possible $Th^n$'s which are determined by them.
\vv \qed 
\sec
%
%
%
%
%
\hh 6. The forcing \par   
\\To contradict the existence of an interpretation we will need generic 
\sc s in every regular cardinal. 
To obtain that we use a simple class forcing. \par
\pd Context. $V\models $ G.C.H \vv
%
%
%
%
\pd \dd 6.1. Let $\la>\ale$  be a regular cardinal
\item{1)} $SC_\la$ := 
$\big\{ f: \ \ f\colon\al\to\{0,1\},\ \al<\la,\ cf(\al)\le\om\  \big\}$ \ 
where each $f$, considered to be a subset of $\al$ (or $\la$), is a \sc. 
The order is inclusion. (So $SC_\la$ adds a generic \sc\ to $\la$). 
\item{2)} $Q_\la$ will be an iteration of the forcing $SC_\la$ with
length $\la^+$ and with support $\le\la$. 
\item{3)} $P$ := $\lan P_\mu$, \qmu : $\mu {\rm \ a\ cardinal >\ale\ } 
\ran$ where \qmu \ is forced to be $Q_\mu$ if $\mu$ is regular,
otherwise it is $\emptyset$. The support of $P$ is sets: each condition in $P$ 
is a function from the class of cardinals to names of conditions where the 
names are non-trivial only for a set of cardinals.
\item{4)} $P_{<\la}, P_{>\la}, P_{\le\la}$ \ are defined naturally. 
For example $P_{<\la}$ is $\lan P_\mu$, \qmu : $\ale<\mu<\la \ran$. \vv
%
%
%
%
\pd Remark 6.2. Note that (if G.C.H holds) $Q_{\la}$ and $P_{\ge\la}$ do not 
add subsets of 
$\la$ with cardinality $<\la$. Hence, $P$ does not collapse cardinals and does 
not change cofinalities, so $V$ and $V^P$ have the same regular cardinals. 
Moreover, for a regular $\la>\ale$ we can split the forcing into 3 parts, 
$P=P_0*P_1*P_2$ where $P_0$ is $P_{<\la}$, $P_1$ is a $P_0$-name of the 
forcing $Q_\la$ and $P_2$ is a $P_0*P_1$-name of the forcing $P_{>\la}$ 
\st $V^P$ and $V^{P_0*P_1}$ have the same $H(\la^+)$.  \par
\\In the next section, when we restrict ourselves to $H(\la^+)$ it will suffice 
to look only in $V^{P_0*P_1}$. \vv
\sec
\hh 7. The contradiction \par
Collecting the results from the previous sections we will reach a 
contradiction from the assumption that there is, in $V^P$, an interpretation 
of $T$ in the monadic theory of a chain $C$. For the moment we will assume 
that the minimal major initial segment $D$ is regular (i.e. isomorphic to 
a regular cardinal), later we will dispose of this by using formal 
shufflings. So we may assume the following:
\pd Assumptions. \par
\\1. $C\in V^P$ interprets $T$ by $\lan d, U_C(\bb X,\bb V), 
E_C(\bb X,\bb Y,\bb V), P(\bb X,\bb Y,\bb V) \ran$. \par
\\2. $D=\la$ is a minimal major initial segment of $C$, $cf(\la)=\la>\om$. \par 
\\3. $\bb R\sb (C-D)$ and \ $S:=\{ \bb A\sb C : \bb A\cap (C-D)=\bb R \}$ \ 
contains an infinite number of nonequivalent representatives of 
$E_C$-equivalence classes.     \par
\\4. There are formulas $U(\bb X,\bb Z), E(\bb X,\bb Y,\bb Z), 
Atom(\bb X,\bb Z), Set(\bb Y,\bb Z)$ and $Code(\bb X,\bb Y,\bb Z)$ in the 
language of the monadic theory of order \st for every $k<\om$ there is a 
sequence $\bb W\sb D$ \st 
$$ I = \lan d, U(\bb X,\bb W), E(\bb X,\bb Y,\bb W), Atom(\bb X,\bb W),
Set(\bb Y,\bb W), Code(\bb X,\bb Y,\bb W) \ran $$
is an interpretation of \tk\ in $D$. \par
\\5. There is an $n<\om$ \st for every $k$ and $\bb W$ as above, 
$Th^n(D,\bb U_{i_1},\bb U_{i_2},\bb U_{i_3},\bb W)$ \ determines the 
truth value of all the interpreting formulas when we replace the variables 
with elements from $\{\bb U_{i_1},\bb U_{i_2},\bb U_{i_3}\}$.  \par 
\\6. $m$ is \st for every $\bb U_{i_1},\bb U_{i_2},\bb U_{i_3}$, \  
from \ $WA^m(D,\bb U_{i_1},\bb U_{i_2},\bb U_{i_3},\bb W)$ 
\ we can compute \ $Th^{n+d}(D,\bb U_{i_1},\bb U_{i_2},\bb U_{i_3},\bb W)$
\ and in particular, the truth value of the interpreting formulas. \par
\\7. Let $N_1$ := $\big| \{ Th^n(C,\bb X,\bb Y,\bb Z) : C$\ is a chain, 
$\bb X,\bb Y,\bb Z \subseteq C\} \big|$. Then (by proposition 3.3 and 
theorem 3.11), for every proper initial segment $D'\sbb D$ there are less 
than $N_1$ $E_C$ -nonequivalent (hence $E$ -nonequivalent) elements, 
coinciding outside $D'$. \par
%
%
%
%
\pd \dd 7.1. The {\sl vicinity} $[\bb X]$ of an element $\bb X$ is the 
collection 
$\{ \bb Y :\ $some element\ $\bb Z\sim\bb Y$\ coincides with\ $\bb X$\  
outside\ some proper (hence minor) initial segment of\ $D\ \}$. \vv
%
%
%
%
\pt Lemma 7.2. Every vicinity $[\bb X]$ is the union of at most $N_1$ 
different equivalence classes.  \vv
\pd Proof. See \gu\ lemma 9.1. \vv \qed
Next we use Ramsey theorem for definining the following functions.
%
%
%
%
\pd Notation 7.3. 
\item{1.} Given  $k<\om$, let $t(k)$ be \st for every sequence
$\bb W\sb D$ of a prefixed length and $a\sb \la$ and for every sequences of 
elements $\lan \bb B_i : i< t(k) \ran$ and $\lan \bb B_s : s\sb t(k) \ran$
there are subsequences $s,s'\sb t(k)$ with $|s'|\ge k$ and $s'\sb s$ \st 
$a-WA^m(D,\bb B_i,\bb B_j,\bb B_s,\bb W)$ is constant for every $i<j\in s'$.  
\item{2.} Given  $k<\om$, let $h(k)$ be \st for every coloring of 
$\big\{(i,j,l) : i<j<l<h(k) \big\}$ into 32 colors, there is a subset $I$ of 
$\{0,1,\ldots,h(k)-1\}$ \st $|I|>k$ and all the triplets 
$\big\{(i,j,l) : i<j<l,\ i,j,l\in I \big\}$ have the same color. \medbreak
%
%
%
%
We are ready now to prove the main theorem:
\pt Theorem 7.4. Assuming the above assumptions we reach a contradiction \vv 
\pd Proof. The proof will be splitted into several steps. 
\medbreak
\\STEP 1: \ Let $K_1 := h(t(3N_1))$ and $K := h(t(2K_1+2N_1))$.\  Let
$\bb R\sb (C-D)$ be \st 
$S := \{ \bb A\sb C : \bb A\cap (C-D)=\bb R \}$ \ contains an infinite number 
of nonequivalent representatives. 
Choose sequences of nonequivalent elements from $S$, 
\ $B := \lan \bb U_i : i<K \ran$, \ and 
$B_1 := \lan \bb V_s : s\sb \{0,1,\ldots,K-1\} \ran$ \ and an appropriate 
$\bb W\sb D$  and interpret $T_K$ on $D$ \st $B$ is the family of ``atoms'' of 
the interpretation and $B_1$ the family of ``sets'' of the interpretation.
\medbreak
\\STEP 2: \ Choose $J:=\{\al_j:j<\la\}\sb\la$ \ 
an $(n+d)$-suitable club witnessing $ATh^{n+d}$ 
for every combination you can think of from the $U_i$'s, the 
$\bb V_s$'s and $\bb W$. \par
\\Now, everything mentioned happens in $H(\la^+)^{V^P}$ 
and, using a previous remark and notations, it is the same thing as 
$H(\la^+)^{V^{P_0*P_1}}$. $P_1$ is an iteration of length $\la^+$ and it 
follows that all the mentioned subsets of $\la$ are added to 
$H(\la^+)^{V^{P_0*P_1}}$ after a proper initial segment of the forcing which 
we denote by $P_0*(P_1|_\be)$. 
So there is a \sc\ $a\sb\la$ in $H(\la^+)^{V^{P_0*P_1}}$ which is added after 
all the mentioned sets, say at stage $\be$ of $P_1$.
\medbreak
\\STEP 3: We will begin now to shuffle the elements with respect to $a$ and 
$J$. \ Let, for $i<j<K$, \  
$k(i,j) := {\rm Min}\{ k : [\bb U_i,\bb U_j]_a^J\sim\bb U_k,\ {\rm or}\ 
k=K \}$. By the definitions of $h$ and $K$ there is a subset 
$s\sb\{0,1,\ldots,K-1\}$ of cardinality at least $K_2:=t(2K_1+2N_1)$ 
\st for every $\bb U_i,\bb U_j,\bb U_l$ with $i<j<l,\ i,j,l\in s$ the
following five statements have the same truth value: \par
\\$k(j,k)=i$, $k(i,k)=j$, $k(i,j)=i$, $k(i,j)=j$, $k(i,j)=k$. Moreover, by 
\gu\ lemma 10.2, if there is a pair $i<j$ in $s$ \st $k(i,j)\in s$ then,  
either for every pair $i<j$ in $s$, $k(i,j)=i$ or for every $i<j$ in $s$, 
$k(i,j)=j$. \medbreak
\\STEP 4: Let $\bb V_s$ be the set that codes 
$\lan \bb U_i : i\in s \ran$. \ By the definitions of $t$ and $K_2$, there 
is a set $s'\sb s$ with at least $K_3:=2K_1+2N_1$ elements and a sequence 
$\lan \bb U_i : i\in s' \ran$ \st for every $r<l$ in $s'$,
$a-WA^m(D,\bb U_r,\bb U_l,\bb V_s,\bb W)$ is constant. \par
\\It follows that for every $r<l$ in $s'$, \ $a-WA^m(D,\bb U_r,\bb V_s,\bb W)$ =  
$a-WA^m(D,\bb U_l,\bb V_s,\bb W)$, and by the preservation theorem 4.12 they 
are equal to $a-WA^m(D,[\bb U_r,\bb U_l]_a^J,\bb V_s,\bb W)$. But  
$\bb V_s$ codes $s$ so $D\models Code(\bb U_r,\bb V_s,\bb W)$, and 
since we can decide from $a-WA^m$ if $Code$ holds, the equality of the 
theories implies that 
$D\models Code([\bb U_r,\bb U_l]_a^J,\bb V_s,\bb W)$.
But by the definition of $Code$ there is $k\in s$ \st 
$[\bb U_r,\bb U_l]_a^J\sim U_k$. So there are $r,l$ in $s$ with $k(r,l)\in s$  
and by step 3 we can conclude that, without loss of generality, for every  
$i<j$ in $s$, \ $[\bb U_i,\bb U_j]_a^J\sim\bb U_i$. \medbreak
\\STEP 5: Note that if $a$ is a \sc\ then $\la\setminus a$ is also a \sc. We will use 
the fact that $a$ is generic with respect to the other sets for finding a 
pair $i<j\in s'$ \st $[\bb U_i,\bb U_j]_{\la\setminus a}^J\sim\bb U_i$ holds as well.
Let $p\in P_0*P_1$ be a condition that forces the value of all the 
theories $a-WA^m(D,\bb U_r,\bb U_l,\bb V_s,\bb W)$ for $r<l\in s'$.
The condition $p$ is a pair $(q,r)$ where $q\in P_0$ and $r$ is a $P_0$-name 
of a function from $\la^+$ to conditions in the forcing $SC_\la$. $r(\be)$ is 
forced by $p$ to be an initial segment of $a$ of height $\ga<\la$ and 
w.l.o.g. we can assume that $\ga=\al_{j+1}\in J$. (So $cf(\ga)=\om$).
As $\ga<\la=D$, $\ga$ is a minor segment. 
Remember that $|s'|\ge K_3=2K_1+2N_1$ and define $s''\sb s'$ 
to be $\big\{ i\in s'$ : $|\{j\in s' : j<i\}|>N_1$, and 
$|\{j\in s' : j>i\}|>N_1 \big\}$. So $|s''|>K_1$. \ Denote by $\fff A\bb B$ 
the element $(\bb A\cap\ga)\cup (\bb B\cap(D-\ga))$. \ We claim that for 
every $i,j,k$ in $s''$, \  
$\bb U_k\sim\ff{[\bb U_i,\bb U_j]_a}\bb U_k$. \par
\\To see that note that by the definition of $s'$ and the preservation 
theorem for $ATh$, \ $p$ \ forces: 
``$Th^{n+d}(D,\ff {[\bb U_i,\bb U_j]_a}\bb U_k,\bb V_S,\bb W)$ = 

\\$Th^{n+d}(D,[\bb U_i,\bb U_j]_a,\bb V_S,\bb W)|_\ga$ +   
$Th^{n+d}(D,\bb U_k,\bb V_S,\bb W)|_{[\ga,\la)}$ = 

{\coo (by $\ga\in J$ and the equality of the a-WA's and the 
                                                       preservation theorem)} 

\\$Th^{n+d}(D,\bb U_i,\bb V_S,\bb W)|_\ga$ + 
$Th^{n+d}(D,\bb U_k,\bb V_S,\bb W)|_{[\ga,\la)}$ = 

{\coo (by $\ga\in J$ and the equality of the ATh's)}

\\$Th^{n+d}(D,\bb U_i,\bb V_S,\bb W)$''. 

\\Hence, since $\bb V_s$ codes $s$, \ 
$\ff {[\bb U_i,\bb U_j]_a}\bb U_k\sim \bb U_l$ for some $l\in s$. \ 
If $l=k$ we are done so assume w.l.o.g that $l<k$. 
Now $U_l\in [U_k]$ and we will show that for every $m<k$, in $s$, 
$U_m\in [U_k]$. Contradiction follows from the choice of $s''$ and lemma  
7.2 (1). 

\\Now $Th^{n+d}(D,\bb U_m,\ff{[\bb U_i,\bb U_j]_a}\bb U_k,\bb W)$ = 

\\$Th^{n+d}(D,\bb U_m,[\bb U_i,\bb U_j]_a,\bb W)|_\ga$ + 
$Th^{n+d}(D,\bb U_m,\bb U_k,\bb W)|_{[\ga,\la)}$. 

\\But $ATh^{n+d}(D,\bb U_m,\bb W)$  = $ATh^{n+d}(D,\bb U_l,\bb W)$. \ 
So there is $\bb Y\sb D$ \st 

\\$Th^n(D,\bb U_m,\bb Y,\bb W)|_\ga$ = 
$Th^n(D,\bb U_l,[\bb U_i,\bb U_j]_a,\bb W)|_\ga$. 

\\We get $Th^n(D,U_m,\fff Y\bb U_k,\bb W)$ = 
$Th^n(D,U_l,\ff{[\bb U_i,\bb U_j]_a}\bb U_k,\bb W)$, \ and since 
$\ff {[\bb U_i,\bb U_j]_a}\bb U_k\sim \bb U_l$ \ the equality of the theories 
implies: \ 
$\fff Y\bb U_k\sim U_m$, \ so $\bb U_m\in[\bb U_k]$.  

\\But by 7.2 (1), \  
$|[U_k]|\le N_1$ and by the choice of $s''$ there are more than $N_1$ 
nonequivalent $U_m$'s with the same property and this is a contradiction. 

So we have proven that  it is possible to replace an initial segment of an 
element with a shuffling of two other elements without changing it's 
equivalence class. (Actually there are $|s''|$  elements like that).
\medbreak
\\STEP 6: \ We are ready to prove that for every $i<j$ in $s''$, \ 
$[\bb U_i,\bb U_j]_a\sim[\bb U_i,\bb U_j]_{\la\setminus a}$. 

\\By step 4 \ $p\fo[\bb U_i,\bb U_j]_a\sim\bb U_i$ \ 
(because it forces equality of theories for a large number of elements).
Remember that $p$ `knows' only an initial segment of $a$, namely only 
$a\cap(j+1)$ where $\ga=\al_{j+1}$.  Since our forcing is homogeneous 
$b$ := $\big(a\cap [0,j+1)\big)\cup\big((\la\setminus a)\cap[j+1,\la)\big)$ \ 
is also generic for all the mentioned sets and parameters, and everything $p$ 
forces for $a$ it forces for $b$. 
\ So $p\fo ``[\bb U_i,\bb U_j]_b\sim \bb U_i''$. 

\\Note that by the preservation theorem 
$Th^n(D,[\bb U_i,\bb U_j]_{\la\setminus a},\bb W)|_\ga$ = 
$Th^n(D,[\bb U_j,\bb U_i]_a,\bb W)|_\ga$ =
$Th^n(D,[\bb U_i,\bb U_j]_a,\bb W)|_\ga$ = 
$Th^n(D,\bb U_i,\bb W)|_\ga$ = $Th^n(D,\bb U_j,\bb W)|_\ga$.

\\It follows that \ 
$Th^n(D,[\bb U_i,\bb U_j]_a,[\bb U_i,\bb U_j]_a,\bb W)|_\ga$ =
$Th^n(D,[\bb U_i,\bb U_j]_{\la\setminus a},[\bb U_i,\bb U_j]_{\la\setminus a},\bb W)|_\ga$.

\\By step 5 (Where we used only the fact that $i,j\in s''$), 
$\ff {[\bb U_i,\bb U_j]_{\la\setminus a}}\bb U_i\sim \bb U_i\sim [\bb U_i,\bb U_j]_b$.  
\ But  
$Th^n(D,\ff {[\bb U_i,\bb U_j]_{\la\setminus a}}\bb U_i,[\bb U_i,\bb U_j]_{\la\setminus a},\bb W)$ =   

\\$Th^n(D,[\bb U_i,\bb U_j]_{\la\setminus a},[\bb U_i,\bb U_j)]_{\la\setminus a},\bb W)|_\ga$ + 
$Th^n(D,\bb U_i,[\bb U_i,\bb U_j]_{\la\setminus a},\bb W)|_{[\ga,\la)}$ =

\\$Th^n(D,[\bb U_i,\bb U_j]_a,[\bb U_i,\bb U_j)]_a,\bb W)|_\ga$ + 
$Th^n(D,\bb U_i,[\bb U_i,\bb U_j]_{\la\setminus a},\bb W)|_{[\ga,\la)}$ =

\\$Th^n(D,\ff {[\bb U_i,\bb U_j]_a}\bb U_i,[\bb U_i,\bb U_j]_b,\bb W)$.

\\But $\ff {[\bb U_i,\bb U_j]_a}\bb U_i\sim \bb U_i\sim [\bb U_i,\bb U_j]_b$. 
\  So it follows by the equality of the theories that
$[\bb U_i,\bb U_j]_{\la\setminus a}\sim[\bb U_i,\bb U_j]_a\sim\bb U_i$ as required.
\medbreak
\\STEP 7: Rename a subsequence of $\lan \bb U_i : i\in s'' \ran$ by 
$\lan \bb A_i : i<2K_1 \ran$ \st for every $i<j<2K_1, r<l<2K_1$ we have:

$(i)$ $a-WA^m(D,\bb A_i,\bb A_j,\bb V_s,\bb W) = 
                               a-WA^m(D,\bb A_r,\bb A_l,\bb V_s,\bb W)$.

$(ii)$ $[\bb A_i,\bb A_j]_a\sim[\bb A_i,\bb A_j]_{\la\setminus a}\sim\bb A_i$.

\\For $i<K_1$ denote by $\bb B_i$ the element that
codes $\bb A_i,\bb A_{2K_1-i-1}$ and look at the sequence 
$\lan \bb B_i : i<K_1 \ran$. \ $K_1$ is large enough so that repeating 
steps 1,2 and 3 we are left with $i<j<K_1$ \st: 

$(iii)$ $a-WA^m(D,\bb A_i,\bb A_{2K_1-i-1},\bb B_i,\bb W) = 
                            a-WA^m(D,\bb A_j,\bb A_{2K_1-j-1},\bb B_j,\bb W)$.

$(iv)$ $[\bb B_i,\bb B_j]_a\sim\bb B_i$ \ or \ $[\bb B_i,\bb B_j]_a\sim\bb B_j$.

\\Now let's shuffle with respect to $a$ and $J$ using clause $(iii)$:

\\$Th^n(D,\bb A_i,\bb A_{2K_1-i-1},\bb B_i,\bb W)$ = 
$Th^n(D,[\bb A_i,\bb A_j]_a,[\bb A_{2K_1-i-1},\bb A_{2K_1-j-1}]_a,
[\bb B_i,\bb B_j]_a,\bb W)$ = 

\\$Th^n(D,[\bb A_i,\bb A_j]_a,[\bb A_{2K_1-j-1},\bb A_{2K_1-i-1}]_{\la\setminus a},
[\bb B_i,\bb B_j]_a,\bb W)$.

\\But $[\bb A_i,\bb A_j]_a\sim\bb A_i$, and by step 6, 
$[\bb A_{2K_1-j-1},\bb A_{2K_1-i-1}]_{\la\setminus a}\sim\bb A_{2K_1-j-1}$ and by 
clause $(iv)$ \ $[\bb B_i,\bb B_j]_a\sim \ \bb B_i$ or $\bb B_j$.

\\So we have, as implied by the equality of $Th^n$ either

\\$\models Code(\bb A_i,\bb B_i,\bb W)\& Code(\bb A_{2K_1-j-1},\bb B_i,\bb W)$ 

\\or 

\\$\models Code(\bb A_i,\bb B_j,\bb W)\& Code(\bb A_{2K_1-j-1},\bb B_j,\bb W)$ 

\\and both cases are impossible!

We have reached a contradiction assuming, in $V^P$, that a well ordered 
chain $C$ interprets $T$ with a minimal major initial segment $D$ which is 
a regular cardinal. 
\vv \qed
\\We still have to prove that there is no interpretation in the case $D$ is 
not a regular cardinal. For that we will use formal shufflings as in 
section 5.
%
%
%
%
\pt Lemma 7.5. The assumption ``$D$ is a regular cardinal'' is not necessary. 
\vv
\pd Proof. Assume first that $D=\de>cf(\de)=\la>\om$. The main point is to 
find 2 elements $\bb A,\bb B$ and a \sc\ $a$ \st 
$[\bb A,\bb B]_a\sim[\bb B,\bb A]_a$ and since $|a|<|A|$, $a$ will be generic  
not with respect to $A,B$ but with respect to sequences of theories of 
length $\la$. We will repeat steps 1 to 7 from the previous proof modifying 
and translating them to the language of formal shufflings. \medbreak
\\STEP 1: We assume $D$ iterprets $T_K$, and choose $\bb W$, K atoms 
$\lan \bb U_i : i<K \ran$ and codings $V_s$ as before.   \medbreak
\\STEP 2: Use notation 2.17*: fix a cofinal sequence in $D$,
$J^* := \lan \be_i : i<\la \ran$, a club $J\sb\la$, 
$J := \lan \al_i : i<\la \ran$\ $(\al_0=\be_0=0)$, and 
$h\colon J^*\rightarrow J$.  W.l.o.g $J$ is an $(n+d)$-suitable club for all 
the combinations of elements we need. (Look at lemma 2.10* and definition 
2.11* for the exact meaning). \par
\\For $\bb k\sb\Big\{\{0,1,\ldots,K-1\}\cup\{(i,j):i<j<K\}\cup
\{s:s\sb \{K-1\}\Big\}$ of length $\le$ 3, let $s_{\bb k}^i$ be the theory 
$Th^{n+d}(\bb U_{\bb k(0)},\ldots,\bb W)|_{[\be_i,\be_{i+1})}$. 
So $Th^{n+d}(\bb U_{\bb k(0)},\ldots,\bb W) =  \sum_{i<\la}s_{\bb k}^i$.

\\Now let ${\cal T}$ denote the set $\{s_{\bb k}^i : \bb k\}\cup
\{\sum_{i<\la}s_{\bb k}^i : \bb k,i\}$.  \ ${\cal T}$ belongs to
$H(\la^+)^{V^{P_0*P_1}} = H(\la^+)^{V^P}$. Call such a ${\cal T}$
{\sl a system of theories}.
\ In $H(\la^+)^{V^P}$ we don't know the $U_i$'s nor the actual ${\cal T}$ but 
we have a set of all the possible systems which must satisfy two sets of 
restrictions: 

\\a) formal restrictions (as in definition 5.1 ).

\\b) material restrictions that reflect the fact that we are dealing with an  
interpretation of $T_K$.  (For example for $\bb k =\lan i,j,\{i,j\}\ran$ the 
theory $\sum_{i<\la}s_{\bb k}^i$ must imply 

\\$Code(\bb X_i,\bb X_{i,j},\bb W)\& Code(\bb X_j,\bb X_{i,j},\bb W)$ ). 

\\So in $H(\la^+)^{V^P}$ we only know that somewhere, (in $H(\de^+)^{V^P}$) 
there are elements that interpret $T_K$ with a system of theories ${\cal T}$.
We scan all the possible systems (they all belong to $H(\la^+)^{V^P}$) and 
show that every one of them leads to a contradiction.

\\Fixing a system ${\cal T}$, let $a\in H(\la^+)^{V^P}$, $a\sb\la$,  
be a generic \sc\ for all the members of ${\cal T}$, which is added at stage 
$\be$ of $P_1$. \medbreak
\\STEPS 3-5: We shuffle the elements with respect to $J$ and $a$ as in 
definition 5.1.(3). The operations are basically the same, but we have to  
translate all the statements to a `formal' language. Just for an example, the 
`formal' meaning of \ $[\bb U_i,\bb U_j]_a^J\sim\bb U_k$ \ is: \ 
``if $s_i = Th^n(\bb U_i,\bb U_k,\bb W)|_{[\be_i,\be{i+1})}$ and 
$t_i = Th^n(\bb U_j,\bb U_k,\bb W)|_{[\be_i,\be{i+1})}$ then 
$Th^n([\bb U_i,\bb U_j]_a^J,\bb U_k,\bb W)$ = $\sum_{i<\la}u_i$ where 
$i\in a\imp u_i=s_i$ and $i\not\in a\imp u_i=t_i$''. 
So $[\bb U_i,\bb U_j]_a^J\sim\bb U_k$ is formally: $\sum_{i<\la}u_i$ implies 
$E(\bb X,\bb Y,\bb W)$.  From this we can easily define 
formally the number $k(i,j)$  as in step 3 in the previous proof. 

\\For choosing a condition $p$ as in step 5, we simply choose a condition 
in $P_0*P_1$ which forces all the `formal' statements we have made. 
This is possible since we are talking about objects of cardinality $\le\la$ 
only. It should be clear that after all the operations we are left with a 
large enough set of elements with some desired properties. Actually if you 
look at the achievements so far, you can note that we didn't use the 
formal theories. $s''$ as in the previous proof can be obtained for any \sc\ 
$a$ so we could have worked in the entire $V^P$ or in $H(\de^+)^{V^P}$. 
But for the next step we need $a$ to be generic.  \medbreak 
\\STEP 6: We have to prove the existence of some $\bb U_i,\bb U_j$ \st $i<j$ 
and $[\bb U_i,\bb U_j]_a\sim[\bb U_j,\bb U_i]_a\sim\bb U_i$. 
Formally we have to prove: 
``if $s_i = Th^n(\bb U_i,\bb U_i,\bb W)|_{[\be_i,\be{i+1})}$ and 
$t_i = Th^n(\bb U_j,\bb U_i,\bb W)|_{[\be_i,\be{i+1})}$ then 
$\sum_{i<\la}u_i$ and $\sum_{i<\la}{u^*}_i$ imply  
$E(\bb X,\bb Y,\bb W)$ where $i\in a\imp u_i=s_i, u^*_i=t_i$ 
and $i\not\in a\imp u_i=t_i, u^*_i=s_i$''.  This follows from the fact that 
$a$ is generic as in step 7 in the previous proof. (Of course, here we can 
not avoid some translation work).   \medbreak
\\STEP 7: We found a \sc\ $a$ and enough elements (at least $K_1$) \st it does 
not matter if we shuffle them with respect to $a$ or with respect to $\la\setminus a$. 
Carry them back to $V^P$ or to $H(\de^+)^{V^P}$ and proceed as before,  
(We don't need the forcing anymore).  

\\The contradiction we have reached proves that ${\cal T}$ can not be realized 
as an interpretation to $T_K$, but since we have chosen it arbitrarily,  
it proves that there is no interpretation at all.  \medbreak
\\STEP 8: We still have to take care of the case ``$D$ is not a well ordered 
chain''.  The only problem is that there may be no first element in $D$, but 
we can fix a $\be_0\in D$ and take into our consideration also theories of 
the form $Th^{n+d}(\bb U_{\bb k(0)},\ldots,\bb W)|_{\be_0}$, but this is  
taken care of in the modified definition of $WA^m$ (look at notation 2.16*). 
Of course all the $K$'s should be computed from the modified definition.
\vv \qed
%
%
%
%
\\Combining 7.4 and 7.5 we get the desired theorem
\pt Theorem 7.6. There is a forcing notion $P$ such that in $V^P$,
Peano arithmetic is not interpretable in the monadic second-order theory of 
chains.
\vv \qed
\sec

\font\ba=cmr8
\font\bs=cmbxti10
\font\bib=cmtt12
\centerline{\bib REFERENCES}  \medbreak
\\ \ {\bf [BaSh]} \ J. B{\ba ALDWIN} and S. S{\ba HELAH}, \ 
{\sl Classification of theories by second order quantifiers}, \ 
{\bs Notre Dame Journal of Formal Logic}, 
\ vol. 26 (1985) pp. 229--303.  \vv
\\ \ {\bf [GMS]} \ Y. G{\ba UREVICH}, M. M{\ba AGIDOR} and S. S{\ba HELAH}, \ 
{\sl The Monadic Theory of $\om_2$}, \ 
{\bs The Journal of Symbolic Logic}, 
\ vol. 48 (1983) pp. 387--398.   \vv
\\ \ {\bf [GU]} \ Y. G{\ba UREVICH},  \ 
{\sl Monadic Second--order Theories}, \ 
{\bs Model Theoretic Logics}, \ 
(J. Barwise and S. Feferman, editors), 
\ Springer--Verlag, Berlin 1985, pp. 479--506   \vv
\\ \ {\bf [GuSh]} \ Y. G{\ba UREVICH} and S. S{\ba HELAH}, \ 
{\sl On the Strength of the Interpretation method}, \ 
{\bs The Journal of Symbolic Logic}, 
\ vol. 54 (1989) pp. 305--323.   \vv
\\ \ {\bf [GuSh1]} \ Y. G{\ba UREVICH} and S. S{\ba HELAH}, \ 
{\sl Monadic Theory of order and topology in ZFC}, \ 
{\bs Ann. Math. Logic}, 
\ vol. 23 (1982) pp. 179--182. \vv
\\ \ {\bf [GuSh2]} \ Y. G{\ba UREVICH} and S. S{\ba HELAH}, \ 
{\sl Interpretating the Second--order logic in the Monadic Theory of Order}, \ 
{\bs The Journal of Symbolic Logic}, 
\ vol. 48 (1983) pp. 816--828.   \vv
\\ \ {\bf [GuSh3]} \ Y. G{\ba UREVICH} and S. S{\ba HELAH}, \ 
{\sl The monadic Theory and the `Next World'}, \ 
{\bs Israel Journal of Mathematics}, 
\ vol. 49 (1984) pp. 55--68.   \vv
\\ \ {\bf [Sh]} \ S. S{\ba HELAH}, \ 
{\sl The monadic Theory of Order}, \ 
{\bs Annals of Mathematics}, 
\ ser. 2, vol. 102 (1975) pp. 379--419.   \vv
\\ \ {\bf [Sh1]} \ S. S{\ba HELAH}, \ 
{\sl Notes on Monadic Logic Part B: Complicatedness for the class of linear orders}, \ 
{\bs Israel Journal of Mathematics}, 
\ vol. 69 (1990) pp. 64--116. \vv

\end